\documentclass[12pt,a4paper]{amsart}
\usepackage{amssymb,graphicx,tikz-cd,mathtools}
\newtheorem{thmA}{Theorem}

\newtheorem{thm}{Theorem}[section]
\newtheorem{cor}[thm]{Corollary}
\newtheorem{lem}[thm]{Lemma}

\theoremstyle{remark}

\def\Z{\mathbb{Z}}
\def\N{\mathbb{N}}
\def\<{\langle}
\def\>{\rangle}
\def\a{\alpha}
\def\b{\beta}

\def\f{F}
\def\T{\mathcal{T}}
\def\G{\Gamma}
\def\onto{\to \hskip -0.35cm \to}
\def\imm{\looparrowright}
\def\Brodsky{Brodski\u\i}
\def\wh{\widehat} 
\def\inac{totally irreducible}

\newenvironment{pfof}[1]{\par\medskip\noindent{\em Proof of #1.}}{\hfill$\square$\par\medskip}

\catcode`\@=11
\def\serieslogo@{\relax}
\def\@setcopyright{\relax}
\catcode`\@=12

\begin{document}

	\begin{abstract}
		We extend several results of Helfer, Wise, Louder and Wilton 
		related to coherence in one-relator groups to the more general setting 
		of one-relator products of locally indicable groups.  
		The methods developed to do so also give rise to a new proof of a theorem of \Brodsky.
	\end{abstract}

\title[coherent 1-relator products]{ 
Coherence  and one-relator products of locally indicable groups}

\author[Howie]{James Howie }
\address{ James Howie\\
Department of Mathematics and Maxwell Institute for Mathematical Sciences\\
Heriot--Watt University\\
Edinburgh EH14 4AS }
\email{ j53howie@gmail.com}


\author[Short]{ Hamish Short } 
\address{ Hamish Short \\
Institut de Mathematiques de Marseille (I2M)\\
39 Rue Joliot--Curie\\
Marseille Cedex 13, France\\
Aix Marseille Univ, CNRS, Centrale Marseille, I2M}
\email{ hamish.short@univ-amu.fr }
 
\thanks{The first named author was supported in part by 
Leverhulme Trust Emeritus Fellowship EM-2018-023$\backslash$9 }
\maketitle 

\section{Introduction} 

In a major recent breakthrough, Louder and Wilton \cite{LW2}  -- and independently Wise \cite{Wi} -- have shown that one-relator groups with torsion are coherent.  In other words, every finitely generated subgroup of such a group has a finite presentation.
This gives a partial answer to an old question of  G. Baumslag \cite{B}.  

A sizeable body of work over the past 40 years, starting with the papers of Brodski\u\i\ \cite{Br,Br0} and the authors \cite{H1,Sh}, has shown that much of one-relator group theory extends to one-relator products of locally indicable groups.  (Recall that a group is {\em locally indicable}  if each of its non-trivial finitely generated subgroups admits an epimorphism onto the infinite cyclic group $\Z$.)  In that spirit, we prove in the current paper the natural analogue of this coherence result, as follows.
 
\begin{thmA}\label{main}
Let $G_\lambda$, $\lambda\in\Lambda$, be a collection of coherent, locally indicable groups, let $S\in *_\lambda G_\lambda$ be a cyclically reduced word of length at least $2$, and let $n>1$ be an integer.  Then the one-relator product
$$G:=\frac{*_\lambda G_\lambda}{\<\<S^n\>\>}$$ 
is coherent.
\end{thmA} 

Baumslag's coherence question remains open in the case of  torsion-free one-relator groups. 
 But the above mentioned theorem in the torsion case is built on earlier work of Helfer and Wise \cite{HW} and of Louder and Wilton \cite{LW}, most of which also applies to torsion-free one-relator groups  and yields partial results in support of the idea that they too are coherent.   We are also able to prove in the current article natural analogues of many of these results in the setting of one-relator products.  We describe these generalizations below. 

We are happy to acknowledge that in the construction of our proofs we have leant heavily on the arguments of Helfer, Wise, Louder and Wilton in the articles cited above, many of which can be readily transported into our framework.  There are of course also some additional difficulties in the more general setting, but we have been able to resolve these.

A two-dimensional CW-complex $Y$ is said to have {\em non-positive immersions} if, for every compact, connected, non-contractible 2-complex $Y'$ admitting an immersion $Y'\imm Y$, the Euler characteristic $\chi(Y')$ is non-positive.

Motivated by Baumslag's conjecture, Helfer and Wise \cite{HW} and independently Louder and Wilton \cite{LW} show that torsion-free one-relator group presentations have the non-positive immersions property, and use this fact in different applications.   For example, Louder and Wilton \cite{LW} show that non-trivial finitely generated subgroups of torsion-free one-relator groups have finitely generated Schur multiplier - indeed the rank of the multiplier is strictly less than the rank of the abelianisation.  As a consequence, one can rule out many incoherent groups such as Thompson's group $F$, the wreath product of $\Z$ with $\Z$, and the direct product of two non-abelian free groups, as subgroups of torsion-free one-relator groups. 

In the present article, 
 we prove  relative analogues of some of these results for one-relator products of locally indicable groups.  (The results also follow for staggered products of locally indicable groups, 
 since these can be constructed as iterative one-relator products.)   
 Now it is easy to show that each component of a 2-complex with non-positive immersions 
 has locally indicable fundamental group  (see \cite[Thm. 1.3]{Wi}), 
 so the local indicability criterion can be omitted from the statement of the  
result on non-positive immersions (Theorem \ref{nonpos} below).

The following construction   -- the reverse of the {\em simple reduction} of \cite{H1} --
occurs throughout the paper, so it is convenient to give it a name.  
Following \cite{HW}, we say that a CW-complex $Y$ is a {\em simple enlargement} of a CW-complex $X$ if
 $Y$ is obtained from $X$ by adjoining a $1$-cell $e$ and at most one $2$-cell $\alpha$, and in the latter case the attaching path $R$ for $\alpha$:
\begin{enumerate}
\item is a closed combinatorial edge-path in $X^{(1)}\cup e$ and involves $e$;
\item is not freely homotopic in $X\cup e$ to a path that crosses $e$ fewer times; and
\item does not represent a proper power in $\pi_1(X\cup e)$.
\end{enumerate}

In practice, we will always assume that $Y$ is connected.  So $X$ has either one or two components -- say $X_1,X_2$  in the latter case.  Then $\pi_1(X\cup e)$ is isomorphic to a free product $\pi_1(X)*\Z$ or $\pi_1(X_1)*\pi_1(X_2)$,  
while if $Y\setminus X$ has a $2$-cell $\alpha$ with attaching map in the homotopy class of $R\in\pi_1(X\cup e)$ then $\pi_1(Y)$ is a {\em one-relator product}
$$\frac{\pi_1(X)*\Z}{\<\<R\>\>}\quad\quad\mathrm{or}\quad\quad\frac{\pi_1(X_1)*\pi_1(X_2)}{\<\<R\>\>}.$$

\begin{thmA}\label{nonpos}
Let $X$ be a $2$-complex  with non-positive immersions, and let $Y$ be a simple enlargement of $X$.
Then $Y$ has non-positive immersions.
\end{thmA}

Applying this to the case where $X$ is one-dimensional, we recover the main result of Helfer and Wise \cite[Theorem 1.3]{HW} and of Louder and Wilton  \cite[Corollary 4]{LW}:

\begin{cor}
Every  torsion-free one-relator group presentation has non-positive immersions.
\end{cor}

The proofs of Theorems \ref{main} and \ref{nonpos} follow a similar pattern to those of Wise in \cite{Wi}, 
of Helfer and Wise in \cite{HW}, and Louder and Wilton in \cite{LW2}.

We also prove an analogue of the theorem of Louder and Wilton about the second Betti number $\beta_2(K)$ of a finitely generated subgroup $K$ of a torsion-free one-relator group.  

Let us say that a group $G$ has the 
{\em second Betti number property} if, for any non-trivial finitely generated subgroup $K$ of $G$, the second Betti number $\beta_2(K)$ of $K$ is strictly less than the first Betti number $\beta_1(K)$. 
Louder and Wilton \cite[Corollary 5]{LW} show that torsion-free one-relator groups have the 
second Betti number property.
Below we prove an analogous result for one-relator products.

\begin{thmA}\label{sb2}
Let $$G:=\frac{\left(\ast_{\lambda\in\Lambda} G_\lambda\right)}{\<\<R\>\>}$$
be a one-relator product of locally indicable groups, each with the 
second Betti number property, where $R\in\ast_\lambda G_\lambda$ is cyclically reduced of length at least $2$, and not a proper power.  Then $G$ has the
second Betti number property.
\end{thmA}

Indeed, we prove a slight generalisation of this theorem as follows.

\medskip
Let $F : \N\to\N$ be a supra-linear function, (i.e. $F(x+y) \ge F(x)+F(y)$ for $x, y \ge 0$; hence in particular
$F(0) = 0$). 
Let $G$ be a finitely presented locally indicable group. We say
that $F$ is a {\em second Betti bounding function} for $G$ if for any non–trivial finitely
generated subgroup $K$ of $G$, we have that 
$$\hat\b(K) := \b_2(K) - \b_1(K) + \b_0(K) \le F(\b_1(K) - \b_0(K)),$$ 
where $\b_i(K)$ denotes the $i$'th Betti number of $K$. 
(As $G$ is locally indicable, $\b_1(K) \ge 1$
when $K$ is finitely generated and non-trivial,
so $\b_1(K)-\b_0(K)$ belongs to the domain $\N$ of $F$ and the above inequality makes sense).

Note that $G$ has the second Betti number property if and only if the zero function is a second Betti bounding function,
and so Theorem \ref{sb2} follows from :  

\begin{thmA}\label{BettiBound}
Let $X$ be a 2-complex with one or two components, each having a locally indicable
fundamental group, and let $Y$ be a simple enlargement of $X$.
If $F : \N\to\N$ is a second Betti bounding function for the fundamental groups
of the components of $X$, then $F$ is a second Betti bounding function for $\pi_1(Y)$.
\end{thmA}

\medskip\noindent{\bf Example.} 
Let $\Lambda$ be a limit group.  Then $\Lambda$ can be constructed from a collection of free abelian groups of finite rank by a series of constructions which are either free products, free products with cyclic amalgamation, or HNN extensions with cyclic amalgamation.  These constructions  are special cases of simple enlargements.
There is a finite upper bound $r>1$ to the rank of abelian subgroups of $\Lambda$.
The map $n\mapsto \frac{n(r-1)(r-2)}2$ is a second Betti bound for $\Z^r$, and hence also for $\Lambda$,  by iterated applications of Theorem \ref{BettiBound}.

In particular, if $\Lambda$ is hyperbolic (or more generally if $\Lambda$ has no abelian subgroup of rank greater than $2$), then  the zero function is a second Betti bound for $\Lambda$, so $\Lambda$ has the second Betti number property.

\bigskip
We also consider the existence or otherwise of certain incoherent subgroups and prove that their non-existence is preserved under simple extensions.  

\begin{thmA}\label{noThompson}
Let $$G:=\frac{\left(\ast_{\lambda\in\Lambda} G_\lambda\right)}{\<\<R\>\>}$$
be a one-relator product of locally indicable groups. Let $K<G$ be a subgroup of $G$ isomorphic to one of the following:
\begin{enumerate}
\item the wreath product $\Z~wr~\Z$;
\item Thompson's group $F$;
\item the direct product of two free groups of rank at least $2$.
\end{enumerate}
Then there exists a unique $\lambda\in\Lambda$ and a unique right coset $G_\lambda g$ of $G_\lambda$ in $G$ such that $K<G_\lambda^g$.
\end{thmA}
 
\noindent{\bf Remarks.}
 Part 3 of Theorem \ref{noThompson} is of particular interest. It complements a result of \Brodsky\ \cite[Theorem 8]{Br} which gives strong restrictions on subgroups of $G$ that decompose as direct products with at least one non-free direct factor.

In the classical case of one-relator groups, Theorem \ref{noThompson} 
says that none of these incoherent groups arise as subgroups of a one-relator group.  
This fact is an easy consequence of the second Betti number property for one-relator groups \cite[Corollary 5]{LW}.

\medskip
The proof of Theorem \ref{noThompson} makes use of the following result of \Brodsky\  \cite[Theorem 6]{Br}. 

\begin{thmA}\label{BrodskiiLemma}  
Let $G_\lambda$, $\lambda\in\Lambda$, be a collection of  locally indicable groups, let $R\in *_\lambda G_\lambda$ be a cyclically reduced word of length at least $2$, and let
$$G:=\frac{*_\lambda G_\lambda}{\<\<R\>\>}.$$
If $g\in G$ and $\lambda,\mu\in\Lambda$ are such that the intersection in $G$ of $G_\lambda$ and $g^{-1}G_\mu g$ is not cyclic, then $\mu=\lambda$ and $g\in G_\lambda$.
\end{thmA}

 It turns out that our methods also yield a new proof of this important result.  Since \cite{Br} is not to our knowledge available online or in translation, we feel that is worthwhile including our proof in this article as well.

The remainder of the paper is organised as follows.  In \S \ref{techstuff} we recall some relevant definitions and previous results,  and note some consequences of them.  In \S \ref{tech} we prove a technical result (Theorem \ref{notree}) which underpins the proofs of our main theorems. 
In \S \ref{tfcase} we prove Theorems \ref{nonpos}, \ref{sb2} and \ref{BettiBound}.  
 In \S \ref{torsioncase} we prove Theorem \ref{main} and note stronger versions of Theorems \ref{noThompson} and \ref{BrodskiiLemma} which hold when the relator is a proper power. 
In \S \ref{Z3} we prove the first 
 two parts of Theorem \ref{noThompson}.  In \S \ref{Blemma} we prove Theorem \ref{BrodskiiLemma}.  
Finally in \S \ref{FxF} 
we prove the last  part of Theorem \ref{noThompson}.

We are grateful to the anonymous referee for 
comments and suggestions which have improved our exposition. 
 
\section{Some technical results we shall need}\label{techstuff}
 
A fundamental result about one relator products of locally indicable groups is the
Freiheitssatz, due independently to the authors and S. \Brodsky.
We shall   use 
this result frequently  -- sometimes explicitly; often implicitly in identifying a free factor with its image in a one-relator product.

\begin{thm}[\cite{Br0},\cite{H3},\cite{Sh}]\label{Freiheitssatz}

Let $A,B$ be locally indicable groups and $R\in A*B$ 
a cyclically reduced word of length at least two. 

The  natural map $A\to \frac{A*B}{\<\< R\>\>}$ is injective.	
	\end{thm}

\bigskip

When  proving Theorems \ref{BrodskiiLemma} and \ref{noThompson}, we shall need the following  decomposition of the
cohomology of a one relator product:

\begin{thm}[\cite{H2}, Theorem 3]\label{cohom}
Let $$G=\frac{*_{\lambda\in\Lambda} G_\lambda}{\<\<R^n\>\>}$$ be a one-relator product of locally indicable groups $G_\lambda$, where $R$ is cyclically reduced of length at least $2$ and not a proper power in $*_\lambda G_\lambda$ and $n\ge 1$.  Let $C$ be the cyclic subgroup of $G$ generated by $R$, and let $M$ be a $\Z G$-module.  Then the restriction maps
$$H^k(G;M) \to H^k(C;M)\times\prod_\lambda H^k(G_\lambda; M)$$ 
are isomorphisms for $k>2$ and an epimorphism for $k=2$. 
\end{thm}

Combining Theorem \ref{cohom} with Shapiro's Lemma (see for instance \cite[III.6.2, page 73]{KB}) we obtain the following:

\begin{cor}\label{Shapiro}
Let  $G$ be as in Theorem \ref{cohom} and let $g\in G$.  If $\mu_1,\mu_2$ are distinct elements of $\Lambda$, 
or if $\mu_1=\mu_2$ and $g\notin G_{\mu_2}$, then the intersection $G_{\mu_1}\cap G_{\mu_2}^g$ 
in $G$ has cohomological dimension at most $1$ (and hence is free).
\end{cor}

\begin{proof}
Let $K:=G_{\mu_1}\cap G_{\mu_2}^g$, let $M$ be a $\Z K$-module, and let $k\ge 2$ be an integer. 
By Shapiro's Lemma we have $$H^k(K;M)\cong H^k(G;Hom_{\Z K} (\Z G,M)).$$
Thus by Theorem \ref{cohom} we have an epimorphism
\begin{eqnarray*}
H^k(K;M) &\onto H^k(C;Hom_{\Z K}(\Z G,M))\times \prod_\lambda H^k(G_\lambda;Hom_{\Z K}(\Z G,M))\\
&\cong \prod_g H^k(K\cap C^g;M) \times \prod_\lambda\prod_g H^k(K\cap G_\lambda^g;M),
\end{eqnarray*}
where $g$ ranges across double-coset representatives for $CgK$, $G_\lambda gK$ respectively.  
But  for the  two  terms on the right hand side of this equation
corresponding to $G_{\mu_1}$ and $G_{\mu_2}^g$,  we have $K=K\cap G_{\mu_1}$ and $K=K\cap G_{\mu_2}^g$.
Hence the diagonal map 
$$H^k(K;M)  \to H^k(K;M)\times H^k(K;M)$$ is also an epimorphism, and it follows that $H^k(K;M)=0$. 
\end{proof}
 
\begin{cor}\label{noZ3}
Let  $G$ be as in Theorem \ref{cohom} and let $K<G$ be a free abelian subgroup of $G$ of rank $r>2$.
Then there exists a unique $\lambda\in\Lambda$ and a unique right coset $G_\lambda g$ of $G_\lambda$ in $G$ such that $K<G_\lambda^g$.
\end{cor}

\begin{proof}
It suffices to prove the case where $r<\infty$.

 We use  Theorem \ref{cohom} and  Shapiro's Lemma as in the proof of Corollary \ref{Shapiro}
 (using $\Z$ coefficients) giving
$$\Z=H^r(K;\Z)\cong H^r(G;\Z G/K)\cong H^r(\<R\>;\Z G/K)\times\prod_\lambda H^r(G_\lambda;\Z G/K)$$
$$\cong \prod_g H^r(\< R\>^g\cap K;\Z)\times \prod_{g,\lambda} H^r(K\cap G_\lambda^g;\Z).$$

Since $K$ has trivial intersection with any conjugate of the finite subgroup $\< R\>$ of $G$, 
it follows that there is precisely one conjugate $G_\lambda^g$ of precisely 
one free factor group $G_\lambda$ such that $H^r(K;\Z)=H^r(K\cap G_\lambda^g;\Z)$.  
In particular $K\cap G_\lambda^g$ has index $1$ in $K$, that is $K<G_\lambda^g$.
\end{proof}

\noindent{\bf Remark.} The above result applies more generally (with the same proof) to the case where $K$ is the fundamental group of any closed, aspherical, orientable manifold of dimension greater than $2$.

\bigskip
A map between CW-complexes is said to be {\em combinatorial} if it maps the interior of each cell homeomorphically onto the interior of a cell (of the same dimension).  An {\em immersion} of CW-complexes is a combinatorial map which is locally injective.

We make frequent use of the following useful fact.  
 (See for example \cite[Lemma 4.1]{LW}, \cite[Lemma 3.1]{H3}  or \cite[Lemma 2.2]{Wi}. )

\begin{lem}\label{fold}
Let $\phi:W\to Y$ be a combinatorial map of CW-complexes, with $W$ compact.   Then $\phi$ factorises as $f\circ\phi'$, where $\phi':W\to Y'$ is surjective and $\pi_1$-surjective, and $f:Y'\imm Y$ is an immersion.
\end{lem}

\medskip
Following Louder and Wilton \cite{LW2}, if $Y=X\cup e\cup\alpha$ is a simple enlargement of a $2$-complex $X$, we define a  
map $f:Y'\to Y$ of $2$-complexes to be a {\em branch map} if it is combinatorial on the complement of $f^{-1}(\alpha)$; 
locally injective in the complement of the preimage  $f^{-1}(\xi)$ of a single point $\xi$ in the interior of $\alpha$; 
and if each 
$2$-cell $\alpha'$ in $f^{-1}(\alpha)$ maps to $\alpha$ as a cyclic branched cover branched over $\xi$.   
We will refer to the degree $n\in\Z_+$ of this branched cover as the {\em branch index} of $\alpha'$.

The following is a natural example of a branch map to $Y$ which forms a key tool in the study of one-relator products in which the relator is a proper power.  For a given positive integer $n$, replace $\alpha$ by a $2$-cell $\alpha_n$ whose attaching path is the $n$'th power of that of $\alpha$.  Let $\widehat{Y}_n$ denote the resulting complex, 
and define $\psi_n:\widehat{Y}_n\to Y$ to be the identity on the complement of $\alpha_n$, 
and on $\alpha_n$ the $n$-fold cyclic cover to $\alpha$, branched over  some point $\xi$.   
Then $\psi_n:\widehat{Y}_n\to Y$ is clearly a branch map.  We call it {\em the   $n$-fold branched cover} of $Y$.

In this paper we will make use of van Kampen diagrams and also their duals, which are known as {\em pictures}, which  were introduced by Rourke in \cite{Rou} and adapted to the relative context by the second author \cite{Sh}.  A picture arises from a continuous map $\Sigma\to X$, where $\Sigma$ is a compact orientable surface and $X$ a $2$-complex, using transversality.  It consists of a finite collection of {\em discs} or {\em (fat) vertices}, whose interiors each map homeomorphically onto the interior of a $2$-cell of $X$, and a properly embedded $1$-submanifold of the complement of the interiors of the discs, each component of which is called an {\em arc}, carries a transverse orientation and is labelled by a $1$-cell of $X$.  A small regular neighbourhood of each arc is mapped to the corresponding $1$-cell in the direction of the transverse orientation.  

If $Y=X\cup e\cup\alpha$ is a simple enlargement, 
then any picture over $Y$ can be made into a {\em relative picture} 
by removing all discs that do not map to $\alpha$, and all arcs that do not map to $e$.

For more details on pictures, and an example of their usefulness in group theory we refer the reader to \cite{DH}, \cite{H4}.
 
\section{The main technical result}\label{tech}
  
Let $X$ be a $2$-complex such that every connected component of $X$ has locally indicable fundamental group. 
Let $Y:=X\cup e\cup\alpha$ be a simple enlargement of $X$, 
and let $R$ denote the closed combinatorial path in $X^{(1)}\cup e$ along which $\alpha$ is attached.

From the definition of simple enlargement, the path $R$ does not represent a proper power in $\pi_1(X\cup e)$.
For each positive integer $n$,
let $\psi_n:\widehat{Y}_n\to Y$ be the $n$-fold branched cover of $Y$, as defined in \S \ref{techstuff}.  In particular $\widehat{Y}_n=X\cup e \cup \alpha_n$ where $\alpha_n$ has attaching map $R^n$.

We subdivide $e$ at  its midpoint $x$, forming two {\em half-edges}, 
and choose $x$ as the basepoint for $Y$.  We orient these half-edges so that $x$ is the initial point of each. 
 The $L$ points where $R$ meets $x$ split it into $L$ closed subpaths $r _1,\dots,r_L$ such that $R=r_1r_2\cdots r_L$.  Let  $g_j:=r_1\cdots r_j$, $j=1,\dots,L$
denote the initial segments of $R$. 
  Weinbaum's Theorem \cite[Corollary 3.4]{H1} 
tells us that no proper closed cyclic subpath of $R$ represents the identity element of $G:=\pi_1(Y)$. 
In particular, if $g_i,g_j$ are two of the initial segments of $R$, with $i\ne j$, 
then $g_i^{-1}g_j=r_{i+1}\cdots r_j$
(if $i<j$) or $g_i^{-1}g_j=r_{i+1}\cdots r_Lr_1\cdots r_j$ (if $j<i$)
is such a cyclic subpath, so $g_i\ne g_j$ in $G$.  Thus the $g_i$  represent pairwise distinct elements of $G$.   

It is known that $G:=\pi_1(Y)=\pi_1(Y,x)$ is locally indicable \cite{H1}, and hence left orderable \cite{BH}.   Suppose that $<$ is a left ordering on $G$.   Then we say that a cyclic subword $V=g_i^{-1}g_j$ of $R$ is {\em minimal} (with respect to $<$) if $V\le U$ for any other cyclic subword $U=g_\ell^{-1}g_k$ of $R$.

  Since a proper cyclic subword  $g_i^{-1}g_j$ of $R$ is non-trivial in $G$, it must be either positive or negative (with respect to any given left ordering). Note that $g_i^{-1}g_j$ is positive if and only if the complementary cyclic subword $g_j^{-1}g_i=(g_i^{-1}g_j)^{-1}$ is negative.

\begin{lem}\label{UV}
 Given any left ordering $<$ on $G$, there is a unique cyclic conjugate of $R$ of the form $U\cdot V$ such that  (with respect to $<$):
\begin{itemize}
\item $V$ is the unique minimal cyclic subword of $R$;
\item each proper initial segment of $U\cdot V$ is positive; and 
\item  each proper terminal segment of $U\cdot V$ is negative.
\end{itemize}
\end{lem}

\begin{proof}

 The finitely many cyclic subwords of $R$ are linearly ordered with respect to $<$,  so  
 we may choose one -- say  $V=g_i^{-1}g_j$ -- which is minimal in $G$ with respect to  $<$.  We claim that the indices
$i$ and $j$ are unique with respect to this property  and hence the minimal cyclic subword $V$ is unique, as required.  
 To see this, let $k\in \{1,\dots,L\}$.  Then by choice of $i,j$ we have
$$g_i^{-1}g_k\ge g_i^{-1}g_j.$$
Left multiplying this inequality by $g_i$ gives $g_k\ge g_j$.  
It follows that $g_j$ is the unique minimal element of $\{g_1,\dots,g_L\}$, 
so the index $j$ is unique, as claimed.  
Hence  also if $k,\ell\in\{1,\dots,L\}$ are such that $g_\ell^{-1}g_k=g_i^{-1}g_j$, 
we have $k=j$ and hence $g_\ell=g_i$, and so also $\ell=i$.  
The index $i$ is therefore also unique.

Clearly $R$ has a cyclic permutation of the form $g_j^{-1}Rg_j=U\cdot V$ where $U$ is also a cyclic subword of $R$.

Any proper initial segment of $U\cdot V$ has the form $g_j^{-1}g_k$ in $G$ for some $k\in\{1,\dots,L\}$, $k\ne j$.  Recall that $g_j$ is the least of the words $\{g_i\}$ with respect to $<$, so we have $g_j<g_k$, and left multiplying by $g_j^{-1}$ gives $1<g_j^{-1}g_k$, as claimed.

Any proper terminal segment of $U\cdot V$ is equal in $G$ to the inverse of a proper initial segment, and so is negative, as claimed.
\end{proof}

Let us assume that $R$ already has the form $U\cdot V$ as  given by Lemma \ref{UV}.
We refer to the edge of $\partial\alpha$ whose midpoint is the starting point of $R$ as the {\em associated $1$-cell} of $\alpha$.

Now suppose that $f:Y'\to Y$ is a branch map.  In particular $f$ is an immersion on the complement of the preimage of a single point in the interior of $\alpha$.  Suppose also that $\alpha'\in f^{-1}(\alpha)$ is a $2$-cell
with branch index $n$. 
Then there are $n$ choices of attaching path for $\alpha'$ that are mapped by $f$ to the path $(UV)^n$.  Each of these paths starts at the midpoint of an edge in $f^{-1}(e)$.  We refer to these edges as the {\em low edges} of $\alpha'$.  Similarly, there are $n$ choices of attaching path that are mapped to $(VU)^n$; each starts at the midpoint of an edge called a {\em high edge} of $\alpha'$.  
 
For each cell $\alpha'\in f^{-1}(\alpha)$ we choose one of the low edges of $\alpha'$ and call it the {\em associated $1$-cell} of $\alpha'$. The attaching path for $\alpha'$  starting at the midpoint of the associated $1$-cell is called  the {\em distinguished attaching path} for $\alpha'$.
Note that, since $f$ is an immersion on a neighbourhood of the $1$-skeleton of $Y'$, no edge can be a low edge (resp. high edge) of more than one $2$-cell in $f^{-1}(\alpha)$.  In particular, distinct $2$-cells in $f^{-1}(\alpha)$ have distinct associated $1$-cells.

	\noindent{\bf Definition.}  We say that an edge $e$ in a combinatorial 2-complex $Y$ is {\it almost-collapsible\/}
if there is a 2-cell $\alpha$ in $Y$ with  attaching map which is  a $k$'th power $S^k$, 
and $e$  appears precisely once in $S$ and does not appear in the attaching path of any $2$-cell other than $\alpha$.
We include the case $k=1$ here, when we sometimes omit the word \lq\lq almost\rq\rq.

If some low (resp. high) edge $e'$ of a $2$-cell $\alpha'\in f^{-1}(\alpha)$ is an almost-collapsible edge 
then  the resulting transformation $Y'\to Y'':= Y'\smallsetminus \{e',\alpha'\}$ is called  a {\em low-edge almost-collapse}  
(resp. a {\em high-edge almost-collapse}). 
  
Note that when $k>1$,  $\pi_1(Y')\cong\pi_1(Y'')*C$, with $C$ cyclic of order $k$.   

\begin{thm}\label{notree}
Let $X$ be a $2$-complex such that every connected component of $X$ has locally indicable fundamental group, 
and let $Y$ be a simple enlargement of $X$.  
Suppose that $f:Y'\to Y$ is a branch map,
with $Y'$ compact and connected, and let 
$A:=Y'\smallsetminus f^{-1}(X)$.
 
Let $Z'$ be the subcomplex obtained from $Y'$ by removing all the $2$-cells in $A$ and their associated $1$-cells.  
Suppose also that, for some component $T$ of $Z'$, $f_*(\pi_1(T))=\{1\}$ in $G=\pi_1(Y)$.
Then $Y'$ can be transformed to $T$ through a sequence of low-edge almost-collapses. 
\end{thm}

\begin{proof}
If $Y=X\cup e$ then $T=Z'=Y'$ and there is nothing to prove.  So for the rest of the proof we consider only the case where $Y=X\cup e\cup\alpha$. 

Let $\T$ denote the collection of subcomplexes of $Y'$ that transform to $T$ through a sequence of low-edge almost-collapses.
Then $T\in\T$ so $\T$ is non-empty.  Clearly $\T$ is partially ordered via inclusion.  Since $f^{-1}(\alpha)$ is finite, it follows that $\T$ must have a maximal element $T'$, say.  The assertion of the theorem is that $T'=Y'$, so we argue by contradiction, beginning from the assumption that $T'\ne Y'$.  Note also that $\pi_1(T')$ is a free product of $\pi_1(T)$ together with a finite number of finite cyclic groups.  Since $f_*(\pi_1(T))=\{1\}$ and $G$ is locally indicable, it follows that $f_*(\pi_1(T'))=\{1\}$ in $G$.

Consider the subset $A'$ of $A$ consisting of those $2$-cells  not in $T'$ 
 whose associated $1$-cells meet $T'$ in either one or both of their endpoints. 
 Let $E$ denote the set of half-edges of associated $1$-cells of $2$-cells in $A'$ having an endpoint in $T'$.  
 Note that the other endpoint of such a half-edge belongs to $f^{-1}(x)$ which is disjoint from $Z'$.  
 We orient each half-edge in $E$ from the endpoint in $f^{-1}(x)$ to the endpoint in $T'$, 
 so that $f$ respects orientation on the half-edges in $E$.

Suppose that $e'\in E$.  Then $e'$ is a half-edge of an associated edge of a $2$-cell $\alpha'\in A'$.  Thus (with a suitable choice of orientation) the distinguished attaching path $R'$  for $\alpha'$ has an initial segment of the form $e'.P.(e'')^{-1}$,
where $P$ is an edge-path in $T'$ and $e''$ is a half-edge of an edge $\widehat{e}$ in $f^{-1}(e)$ which is not contained in $T'$.  Now $\widehat{e}$
is the associated $1$-cell of a $2$-cell in $A'$, since otherwise it is contained in $Z'$ and hence in $T'$.
Since $e''$ has an end-point in $T'$ we have $e''\in E$.  Note that $e''$ is uniquely determined by $e'$. We call $e''$ the
{\em successor} of $e'$ and write $e''=\sigma(e')$.  Thus $\sigma:E\to E$ is a well-defined map.

Here are some remarks about this map $\sigma$.

\begin{enumerate}
\item\label{one} $\sigma(e')\ne e'$.  For suppose that $Q:=e'.P.(e')^{-1}$ is a subpath of the attaching path $R'$
of a $2$-cell $\alpha'\in f^{-1}(\alpha)$.  
Since $R'$ is cyclically reduced, it follows that $P$ is non-empty and $Q$ is not the whole of $R'$. 
 Moreover $f_*(Q)$ is an initial segment of $R^{\pm n}$ but is not a power of $R$, and so $f_*(Q)>1$ in $G$,
 by Lemma \ref{UV}.
However, $P$ is a closed path in $T$ and by hypothesis $f_*(\pi_1(T))=\{1\}$ in $\pi_1(Y)$.  
Hence also
$ f_*(Q)=f_*(e'.P.(e')^{-1})=1$ in $G$.  This gives a contradiction.
 
\item\label{three} Suppose that $e'\in E$, and that $e'$ and $e''=\sigma(e')$ are not half-edges of the same edge.  Then
the initial segment $Q=e'.P.(e'')^{-1}$ of $R'=\partial\alpha'$ in the definition is proper and non-empty. 
If $f(Q)=R^k$ for some $k$ then $0<k<n$ and so $e'' $ is a half-edge of another low edge of $\alpha'$. Since it is also the associated $1$-cell of a $2$-cell $\alpha''\in f^{-1}(\alpha)$, and since it cannot be a low edge of two distinct $2$-cells, we must have $\alpha''=\alpha'$ and $k\equiv 0 ~\mathrm{mod}~n$, contrary to assumption.
Hence $f(Q)$ is not a power of $R$.
As in Remark \ref{one} above 
we have $ f_*(Q)>1$ in $G$.

\item\label{four} Since $E$ is finite, any chain of the form $e_1, e_2=\sigma(e_1), e_3=\sigma(e_2),\dots$ in $E$ must contain a loop.
Without loss of generality let us suppose that $e_n=e_1$ for some $n$.  
If, for each pair $e_j,e_{j+1}$ in this sequence, 
$e_j$ and $e_{j+1}$ are half-edges of associated $1$-cells of {\em distinct} cells of $f^{-1}(\alpha)$,
then by Remark \ref{three} we have paths $Q_j=e_j.P_j.e_{j+1}^{-1}$ where $P_j$ is an edge--path in $T'$ 
and $f(Q_j)$ is an initial segment of  $(UV)^{\pm n}$ 
-- and so $f_*(Q_j)>1$ in $G$  by Lemma \ref{UV}.  
Hence $f_*(Q)>1$ where $Q=Q_1.Q_2.\cdots .Q_{n-1}$.  On the other hand
$f_*(Q)=f_*(e_1.P.e_1^{-1})$  where $P=P_1.P_2.\cdots .P_{n-1}$ is a closed path in $T'$ 
and $f_*(\pi_1(T'))=\{1\}$ in $G$.
Hence $f_*(Q)=1$ in $G$.  
This contradiction shows that, if  $E$ is not empty, there must be a 
$2$-cell $\alpha'\in A'$ and a half-edge $e'$ of the associated $1$-cell of $\alpha'$, such that $\sigma(e')$ is also a half-edge of the associated $1$-cell of $\alpha'$.

\item\label{two} By assumption, $T'\ne Y'$ and so $E$ is non-empty.  By Remark \ref{four} there exist $\alpha'\in A'$ and $e'\in E $ such that $e'$ and $e'':=\sigma(e')$ are half-edges of the associated $1$-cell $\widehat{e}$ of $\alpha'$.  By Remark \ref{one} we cannot have $\sigma(e')=e'$, so $e'$ and $e''$ are the two distinct  half-edges of $\widehat{e}$.   
In this case the attaching path $R'$ of $\alpha'$ is a power of $Q=e'.P.(e'')^{-1}$.  So $f_*(Q)=1$ 
in $G$.  Hence $\widehat{e}$ is an almost-collapsible face of $\alpha'$, and $T'\cup\widehat{e}\cup\alpha'$ is a subcomplex of
$Y'$ which admits a low-edge almost-collapse to $T'$.

\end{enumerate}
This gives us the desired contradiction to the maximality of $T'$ in $\T$ and
 completes the proof. 
\end{proof}
 
In the torsion-free case, we have the following stronger version of Theorem \ref{notree}, which we will apply in \S \ref{tfcase} below.

\begin{thm}\label{notreetf}
Let $X$ be a $2$-complex such that every connected component of $X$ has locally indicable fundamental group, and let $Y$ be a simple enlargement of $X$. 
Suppose that $f:Y'\imm Y$ is an immersion,
with $Y'$ compact and connected, and let  
 $A:=Y'\smallsetminus f^{-1}(X)$.
Let $Z'$ be the subcomplex obtained from $Y'$ by removing all the $2$-cells in $A$ and their associated $1$-cells.  
Suppose also that, for some component $T$ of $Z'$, $f_*(\pi_1(T))=\{1\}$ in $G$.

Then $Y'$ collapses to $T$ through a sequence of low-edge collapses. 
\end{thm}

\begin{proof}
The map $f$ is an immersion, hence a branch map with no branch points.  
By Theorem \ref{notree}, $Y'$ can be transformed to $T$ through a sequence of low-edge almost-collapses.   
But since there are no branch points of $f$ in $Y'$, each low-edge almost-collapse is a genuine collapse.  
The result follows.
\end{proof}

\section{The torsion-free case and the Betti number property}\label{tfcase}

In this section we consider one-relator products of locally indicable groups in which the relator is not a proper power.
 We prove Theorems \ref{nonpos} and \ref{BettiBound}, and then show how Theorem \ref{sb2} follows from Theorem \ref{BettiBound}.   We end the section with a slightly strengthened version of the Freiheitssatz (Theorem \ref{Freiheitssatz}), and an application to complexes which immerse into simple enlargements.

\begin{pfof}{Theorem \ref{nonpos}}

\noindent{\bf Case 1:} $Y=X\cup e$ where $e$ is a $1$-cell.

Let $f:Y'\imm Y$ be an immersion with $Y'$ compact and connected, and suppose that $\chi(Y')>0$.  Let $X':=f^{-1}(X)$.  Since $X$ has
the non-positive immersions property, each component of $X'$ either is contractible or has non-positive Euler characteristic.
In particular each component of $X'$ has Euler characteristic at most $1$.  Since $Y'$ is connected and constructed from $X'$
by attaching $1$-cells, it follows that $Y'$ has Euler characteristic at most $1$, and hence by hypothesis $\chi(Y')=1$.
But $\chi(Y')=1$ implies that each component of $X'$ has Euler characteristic $1$ (and hence is contractible) and that these components are connected in a tree-like manner by the edges in $f^{-1}(e)$ to form $Y'$.  Hence $Y'$ is also contractible, as claimed.

\medskip\noindent
{\bf Case 2:} $Y=X\cup e\cup\alpha$ where $e$ is a $1$-cell and $\alpha$ is a $2$-cell.

Let $f:Y'\imm Y$ be an immersion with $Y'$ compact and connected, and suppose that $\chi(Y')>0$.  
Let $Z'$ be as in Theorem \ref{notreetf}.  
Since $Z'$ is obtained from $Y'$ by removing equal numbers of $1$- and $2$-cells, we have $\chi(Z')=\chi(Y')>0$, and hence there is a component $T$ of $Z'$ with $\chi(T)>0$.
   But $f(T)\subset X\cup e$ and $X\cup e$ has the non-positive immersions property by Case 1 above,  so $T$ is contractible.  By Theorem \ref{notreetf}, $Y'$  collapses to $T$, so $Y'$ is also contractible, as claimed.
\end{pfof}

\begin{pfof}{Theorem \ref{BettiBound}}
Let $K$ be a finitely generated subgroup of $\pi_1(Y)$.  Suppose that $K$ has first Betti number $m$ and second Betti number $\beta_2(K)\ge n:=m+F(m-1)$.  Suppose that $K$ can be generated by $d$ elements (where of course $d\ge m$). 
We will show that $K=1$, from which the assertion of the Theorem follows.

  We can construct an epimorphism $\widehat{K}\to K$ for some group $\widehat{K}$ 
  which has first Betti number $m$ and possesses a $d$-generator, $(d-m)$-relator presentation $\mathcal{P}$.  
  Let $F$ be the free group on these $d$ generators, and $\mathcal R$ the kernel of the epimorphism $F\twoheadrightarrow\widehat{K}\twoheadrightarrow K$.  
  Then we can choose relations $r_1,\dots,r_n\in \mathcal{R}\cap [F,F]$ 
  which are linearly independent in $$\frac{\mathcal{R}\cap [F,F]}{[\mathcal{R},F]}\cong H_2(K),$$ 
  and add them to the presentation $\mathcal{P}$ to get a new presentation $\mathcal{P}'$ of a group $K'$, 
  also admitting an epimorphism to $K$.

This epimorphism can be realised by a combinatorial map $f:V\to Y$ for some subdivision $V$ of $\mathcal{P}'$.
Lemma \ref{fold}
gives a factorization of $f$ as $g\circ \bar{f}$, where $g:Y'\imm Y$ is an immersion, 
and $\bar{f}:V\to Y'$ is a surjective combinatorial map which is also surjective on $\pi_1$,
and hence also on $H_1$.
 Hence $Y'$ is connected with first Betti number $\b_1(Y')\le m$.   On the other hand, since $g_*(\pi_1(Y'))=f_*(\pi_1(V))=K$ has first Betti number $m$, we also have $\b_1(Y')\ge m$ and hence $\b_1(Y')=m$.
Now the map $\mathbb{Z}^n\cong H_2(V)\to H_2(K')\to H_2(K)$ is injective and factors through 
$H_2(Y')$.  Hence $Y'$ has second Betti number $\ge n$.

Let $Z'$ be the subcomplex of $Y'$ defined in Theorem \ref{notreetf}.  Since $1\ne K=f_*(\pi_1(V))=g_*(\pi_1(Y'))$, it follows from Theorem \ref{notreetf} that $g_*(\pi_1(T))\ne 1$ in $\pi_1(Y)$ for each component $T$ of $Z'$.

For each component $X_j$ of $X':=g^{-1}(X)$, embed $X_j$ (for example via a mapping cylinder construction) into a
classifying space $\overline{X_j}$ for the subgroup $K_j:=g_*(\pi_1(X_j))$ of $K$, in such a way that the embedding map realises the given map $g_*$ on fundamental groups.  Let $\overline{X}$ denote the disjoint union of these $\overline{X_j}$ for all the components $X_j$ of $X'$.  Form $\overline{Z}$ and $\overline{Y}$ from $Z',Y'$ respectively by adjoining $\overline{X}$ along $X'$.  We can extend the map $g:Y'\to Y$ to a map from $\overline{Y}$ to a classifying space $\widehat{Y}$ for $\pi_1(Y)$, in such a way that the restriction to each $\overline{X_j}$ factors through the covering of $\widehat{Y}$ corresponding to $K_j$.  By hypothesis, each non-trivial $K_j$ has the property that 
$$\b_2(K_j)-\b_1(K_j)+1\le F(\b_1(K_j)-1).$$
Since $\overline{X_j}$ is a classifying space for $K_j$, it follows that each non-contractible $\overline{X_j}$ satisfies
$$\beta_2(\overline{X_j})-\b_1(\overline{X_j})+1\le F(\beta_1(\overline{X_j})-1).$$
Without loss of generality, suppose that there are $J$ non-contractible components $\overline{X_j}$, $1\le j\le J$ of $\overline{X}$, and $C$ contractible components $\overline{X_c}$, $J+1\le c\le J+C$.  Then
\begin{eqnarray*}
\beta_2(\overline{X}) - \beta_1(\overline{X})+\beta_0(\overline{X}) &= \sum_{j=1}^{J+C} \beta_2(\overline{X_j}) - \beta_1(\overline{X_j})+1\\
&=C + \sum_{j=1}^J \beta_2(\overline{X_j}) - \beta_1(\overline{X_j})+1\\
&\le C+ \sum_{j=1}^J F(\beta_1(\overline{X_j})-1).
\end{eqnarray*}
As $F$ is supra-linear, it follows that
\begin{equation}\label{star}
\beta_2(\overline{X}) - \beta_1(\overline{X})+\beta_0(\overline{X}) \le C+F\left(\sum_{j=1}^J \beta_1(\overline{X_j})-1\right).
\end{equation}

Now $\overline{Z}$ is constructed from $\overline{X}$ by adding a finite number -- say $\ell$ -- of $1$-cells.
Then $\overline{Y}$ is constructed from $\overline{Z}$ by adding equal numbers of $1$- and $2$-cells.
It follows that 
\begin{equation}\label{eqtwo}
\beta_2(\overline{Y})-\beta_1(\overline{Y})+\beta_0(\overline{Y})=\beta_2(\overline{Z})-\beta_1(\overline{Z})+\beta_0(\overline{Z})=\beta_2(\overline{X})-\beta_1(\overline{X})+\beta_0(\overline{X})-\ell.
\end{equation}

Since $\overline{Y}$ is constructed from $\overline{X}$ by adding $1$-cells and $2$-cells, we have $H_3(\overline{Y},\overline{X})=0$, and so $H_2(\overline{X})\to H_2(\overline{Y})$ is injective, from the long exact homology sequence.  In particular $\beta_2(\overline{Y})\ge \beta_2(\overline{X})$.

From (\ref{eqtwo}), recalling that $\overline{Y}$ is connected, we have 
$$\b_1(\overline{Y})-1  = \b_2(\overline Y) - \b_2(\overline X) + \b_1(\overline{X})-\b_0(\overline{X})+\ell\ge \b_1(\overline{X})-\b_0(\overline{X})+\ell$$
and as $F$ is non-decreasing we obtain
$$F(\b_1(\overline{Y})-1)\ge F(\b_1(\overline{X})-\b_0(\overline{X})+\ell).$$

We assumed that $F$ is not a second Betti bound for $\pi_1(\overline{Y})$, so from (\ref{star}) and (\ref{eqtwo})
we derive
\begin{align*}
\ \ \ \ F(\b_1(\overline{Y})-1)
&< \b_2(\overline{Y})-\b_1(\overline{Y})+\b_0(\overline{Y})\\
&=\b_2(\overline{X})-\b_1(\overline{X})+\b_0(\overline{X})-\ell\\
&\le C-\ell+F\left(\sum_{j=1}^J (\b_1(\overline{X_j})-1)\right)\\
&=(C-\ell)+F(\b_1(\overline{X})-\b_0(\overline{X})+C)\\
&=(C-\ell)+F(\b_1(\overline{X})-\b_0(\overline{X})+\ell+(C-\ell)).
\end{align*}

But $F$ is non-decreasing and 
$$\b_1(\overline{Y})-1\ge \b_1(\overline{X})-\b_0(\overline{X})+\ell\ge\b_1(\overline{X})-\b_0(\overline{X})+\ell+(C-\ell))$$
if $C-\ell\le 0$, giving a contradiction.

It follows that $C>\ell$.

Thus there is a component $\overline{T}$ of $\overline{Z}$ such that the number $k_T$ of contractible components of $\overline{X}$ contained in $\overline{T}$ is strictly greater than the number $\ell_T$ of $1$-cells of $\overline{Z}\setminus\overline{X}$ contained in $\overline{T}$.  But $\overline{T}$ is connected, and it follows that $k_T=\ell_T+1$, and that  $\overline{T}$ is also contractible.  Hence $T:=\overline{T}\cap Z'$ is a component of $Z'$ such that $g_*(\pi_1(T))=\{1\}$ in $\pi_1(Y)$.
By Theorem \ref{notreetf}, $Y'$ collapses to $T$, so $K=g_*(\pi_1(Y'))=1$ in $\pi_1(Y)$.
This completes the proof.
\end{pfof}
 
\begin{pfof}{Theorem \ref{sb2}}
We first reduce to the two-factor case as follows.  
Choose a partition $\Lambda=\Lambda(1)\sqcup\Lambda(2)$ such that $R$ contains a letter from $G_\lambda$ for at least one $\lambda\in\Lambda(1)$, and a letter from $G_\mu$ for at least one $\mu\in\Lambda(2)$.  
Then $G$ is a one-relator product of $A:=\ast_{\lambda\in\Lambda(1)}G_\lambda$ and $B:=\ast_{\mu\in\Lambda(2)} G_\mu$.
Let $F:\N\to\N$ be the zero function: $F(n)=0~\forall~n$.  
Then, as pointed out in the Introduction, a group $H$ has the second Betti number property if and only if $F$ is a second Betti bounding function for $H$.  In particular, $F$ is a second Betti bounding function for each $G_\lambda$ -- and hence also for $A$ and $B$ by the homological properties of free products.
By Theorem \ref{BettiBound} it follows that $F$ is also a second Betti bounding function for $G$, and hence $G$ has the second Betti number property, as required.
\end{pfof}

\smallskip
We end this section with a slightly strengthened version of the Freiheitssatz  (Theorem \ref{Freiheitssatz}), 
and an application to complexes which immerse into simple enlargements.

\begin{lem}\label{StrongFHS}
Let $Y:=X\cup e \cup\alpha$ be a simple enlargement of a $2$-complex $X$, 
 where every fundamental group of $X$ is locally indicable.  Let  $\phi:\Delta\to Y$ be a reduced van Kampen diagram
which is either spherical or a disk diagram.  If $\phi^{-1}(\alpha)\ne\emptyset$ then $\Delta$ is a disk diagram, and there is a $2$-cell $\alpha'$ mapping to $\alpha$ such that the low-edge of $\alpha'$ lies on the boundary of $\Delta$.
\end{lem}

\begin{proof}
Suppose not.  Then there is a sequence of $2$-cells $\alpha_1,\dots,\alpha_n\in\phi^{-1}(\alpha)$ such that for each $j$ the low edge $e_j$  of $\alpha_j$ is also on the boundary of $\alpha_{j+1}$ (indices modulo $n$).  
By Lemma \ref{UV} and the definition of low edge, either of the two paths $P_{j+1}$ in $\partial\alpha_{j+1}$ from the midpoint of $e_j$ to that of $e_{j+1}$ represents an element $U_{j+1}=[\phi(P_{j+1})]<1$ in $\pi_1(Y)$ (with respect to a fixed left ordering of $\pi_1(Y)$).  Hence their product $$U:=U_1\cdot U_2\cdots U_n<1$$  in $\pi_1(Y)$.

But this contradicts the fact that $U$ is represented by $\phi(P)$ where $$P:=P_1\cdot P_2\cdots P_n$$ is a closed path in the simply-connected space $\Delta$.
\end{proof}

\begin{cor}\label{pi1-inj}
Let $Y$ be as in Lemma \ref{StrongFHS}, and let $f:Y'\imm Y$ be an immersion, with $Y'$ compact.   Let $Z'$ denote the subcomplex of $Y'$ defined in Theorem \ref{notreetf}, and let $\{Z_j\}_{j\in J}$ be its components.  Then 
\begin{enumerate}
\item The inclusion-induced map $\pi_1(Z_j)\to\pi_1(Y')$ is injective for each $j\in J$; and
\item $\pi_2(Y')$ is generated as a $\Z\pi_1(Y')$-module by the images of $\pi_2(Z_j)$ for  all $j\in J$.
\end{enumerate}
\end{cor}

\begin{proof}
\begin{enumerate}
\item
Apply Lemma \ref{StrongFHS} to  a reduced van Kampen diagram in the plane, with boundary labelled by a path in $Z_j$ representing an element $w\in\mathrm{Ker}(\pi_1(Z_j)\to\pi_1(Y'))$.  Since the low edges of $2$-cells in $f^{-1}(\alpha)$ are by definition excluded from $Z'$, the Lemma says that there are no $2$-cells of $f^{-1}(\alpha)$ in $\Delta$.  Hence $\phi(\Delta)\subset Y'\setminus f^{-1}(\alpha)$.  Since $Y'\setminus f^{-1}(\alpha)$ is formed by adding $1$-cells to $Z'$, the map $\pi_1(Z_j)\to\pi_1(Y'\setminus f^{-1}(\alpha))$ is injective, and so the boundary label $w$ of $\Delta$ is already trivial in $\pi_1(Z_j)$.
\item
  Applying Lemma \ref{StrongFHS} to the case where $\Delta$ is a spherical diagram, we see that no $2$-cell of $\Delta$ maps to $\alpha$ under $f\circ\phi$.  It follows that $\pi_2(Y')$ is generated as a $\Z\pi_1(Y')$-module by the image of $\pi_2(Y'\setminus f^{-1}(\alpha))$. Since $Y'\setminus f^{-1}(\alpha)$ is formed by adding $1$-cells to $Z'$, we deduce that in fact $\pi_2(Y')$ is generated as a $\Z\pi_1(Y')$-module by the images of $\pi_2(Z_j)$ for all $j\in J$.
\end{enumerate}
\end{proof}

\section{The torsion case}\label{torsioncase}
 
In this section we prove Theorem \ref{main}, and also stronger versions of Theorems \ref{noThompson} and \ref{BrodskiiLemma} in the torsion case.  The proof of Theorem \ref{main} requires a number of preliminary results.

\begin{thm}\label{p-power}
Let $f:Y'\to Y$ be as in Theorem \ref{notree}, let $p$ be a prime and let $F$ be the field of order $p$. Then the number of $2$-cells of $f^{-1}(\alpha)$ which are attached along $p$'th powers is bounded above by the $F$-dimension of $H_1(Y',F)$.
\end{thm}

\begin{proof}
As in Theorem \ref{notree}, let $Z'$ be formed from $Y'$ by deleting all the $2$-cells of $f^{-1}(\alpha)$ together with their associated $1$-cells.  Suppose first that some component $T$ of $Z'$ has zero first Betti number.   Then $f_*(\pi_1(T))=\{1\}$ in $G$ since $G$ is locally indicable.   By Theorem \ref{notree}, there is a sequence of low-edge almost-collapses transforming $Y'$ to $T$.  Necessarily this sequence of almost-collapses involves every $2$-cell in $f^{-1}(\alpha)$.
Suppose that $f^{-1}(\alpha)=\{\alpha'_1,\alpha'_2,\dots,\alpha'_N\}$ and that $\alpha'_j$ is attached along a $q(j)$'th power for each $j$ (where $q(j)\ge 1$).   Suppose also that $p|q(j)$ for $j\le J$.  Then  $\pi_1(Y')$ is the free product of $\pi_1(T)$ and the cyclic groups $C(j)$ of order $q(j)$ for $j=1,\dots,N$.   Hence $H_1(Y',F)$ is the direct sum of $H_1(T,F)$ and $H_1(C(j),F)$ for $j=1,\dots,N$.  Since $p|q(j)$ for $j\le J$, at least $J$ of these direct factors are isomorphic to $F$, and the result follows in this case.

\medskip
Hence we may assume that every component of $Z'$ has positive first Betti number.   Suppose that there are $K$ components in $Z'$, and $N$ $2$-cells in $f^{-1}(\alpha)$, of which $J$ are attached along $p$'th powers.  Let $Z'':=Z'\cup Y^{(1)}$.   Then $Z''$ is obtained from $Z'$ by adding $N$ $1$-cells, of which $K-1$ are required to make $Z''$ connected, and the remaining $N-K+1$ contribute to the first Betti number.   So $Z''$ has first Betti number at least $K+(N-K+1)=N+1$.  Hence also $H_1(Z'',F)$ has $F$-dimension at least $N+1$.  Finally, $H_1(Y',F)$ is the quotient of $H_1(Z'',F)$ by the subspace $V$ spanned by the images of the attaching paths for the $2$-cells in $f^{-1}(\alpha)$.  Those attaching paths which are $p$'th powers have image $0$ in $H_1(Z'',F)$, so $V$ has dimension at most $N-J$.  Hence $H_1(Y',F)$ has dimension at least $(N+1)-(N-J)=J+1$, and the result follows.
\end{proof}
 
We next consider immersions $f:Y'\imm \widehat{Y}_n$ where $n>1$.   
The $2$-cell $\alpha_n$ of $\widehat{Y}_n$ has attaching map $R^n$, where $R$ is not a proper power.  
It follows that each $2$-cell $\alpha'\in f^{-1}(\alpha_n)$ has an attaching map of the form $S^p$, where $S$ is not a proper power, $f(S)=R^q$, and $pq=n$.

\begin{thm}\label{5beta}
Let $n>1$ and suppose that $f:Y'\imm \widehat{Y}_n$ is an immersion, where $Y'$ is compact and connected, with first Betti number $\b$.  Suppose that none of the $2$-cells in $f^{-1}(\alpha_n)$ is attached by an $n$'th power.  
If $Y'$ has no almost-collapsible edges, then the number of $2$-cells in $f^{-1}(\alpha_n)$ is bounded above by $5\b$.
\end{thm}

\begin{proof}
The composite $\psi_n\circ f:Y'\to Y$ is a branch map.  Let $Z''$ be the subcomplex of $Y'$ obtained by removing all the $2$-cells in $f^{-1}(\alpha_n)=(\psi_n\circ f)^{-1}(\alpha)$, together with all their low edges.  Note that there are at least $2$ distinct low edges for each $2$-cell in $f^{-1}(\alpha_n)$, since its attaching map is not an $n$'th power.

Now let $C$ be a component of $Z''$ with first Betti number $0$.   (We call such a component a {\em treeoid}.)  Then $(\psi_n\circ f)_*(\pi_1(C))$ is a finitely generated subgroup of $G$ with first Betti number $0$, and hence trivial.
If $C=Y'$ then $f(Y')=f(C)\subset X\cup e$, and there are no $2$-cells in $f^{-1}(\alpha_n)$.  
So we may assume that $C\ne Y'$.  
Since $Y'$ is connected, there must be one or more low-edges of $2$-cells in $f^{-1}(\alpha_n)$ that meet $C$.

We claim that $C$ is incident to at least $5$ half-edges of low edges of $2$-cells of $f^{-1}(\alpha_n)$. 
  
Arguing as in the proof of Theorem \ref{notree}, 
for each such half-edge $e'$ there is another half-edge $\sigma(e')$ 
and a cyclic subpath of the attaching path of some $\alpha'\in f^{-1}(\alpha_n)$ 
of the form $Q:=e'\cdot P\cdot \sigma(e')^{-1}$ with $P$ a path in $C$.   
Then $f(Q)$ is an initial segment of the attaching path $R^n=(UV)^n$ of $\alpha_n$, 
so  by Lemma \ref{UV} $f_*(Q)\ge 1$ in the left ordering of $G$, with equality 
if and only if $f(Q)$ is a power of $R$ - necessarily $R^{\pm 1}$ since $P$ cannot contain a low edge of $\alpha'$.

Still arguing as in the proof of Theorem \ref{notree}, 
any chain of half-edges $e_1, e_2, \dots$ of low edges of cells in $f^{-1}(\alpha_n)$ 
that are incident at $C$ with $\sigma(e_i)=e_{i+1}$ must contain a loop, 
and the only possibility for such a loop is a pair $e_1,e_2$ with $e_1=\sigma(e_2)$, $e_2=\sigma(e_1)$, and $f(Q)=R^{\pm 1}$. 

 Now the path $Q=e_1\cdot P\cdot e_2^{-1}$ must contain a high edge $\hat{e}$ of $\alpha'$.  
 Since $\hat{e}$ is not an almost-collapsible edge of $\alpha'$, 
 there must be another subpath $P'$ in $C$ of the attaching path of a $2$-cell 
 $\alpha''\in f^{-1}(\alpha_n)$ that passes through $\hat{e}$.   
 Note that $P'$ is a cyclic 
 subword of $R^{\pm 1}$.
 
 It is not possible for $P'$ to begin  at the edge  $e_1$, 
 else there would be a loop formed from an initial segment $P'_1$ of $P'$ to the edge $\hat e$ together with
 the initial segment $P_1$ of $P$ from $e_1$ to the edge $\hat e$.
Since $f_*(C)=1$, it follows that $f_*(P'_1)=f_*(P_1)=U=V^{-1}$.  
	By uniqueness of low and high edges, it then follows that $e_1$ and $\hat e$ are low and high edges respectively of $\alpha''$, 
and hence that $\alpha''=\alpha'$ and so $P'=P$ and the edge $\hat e$ would be almost-collapsible, a contradiction.

 Essentially the same argument shows that  $P'$ does not begin or end at either of the edges $e_1,e_2$.

\begin{figure}
	\centerline{\includegraphics[scale=0.3]{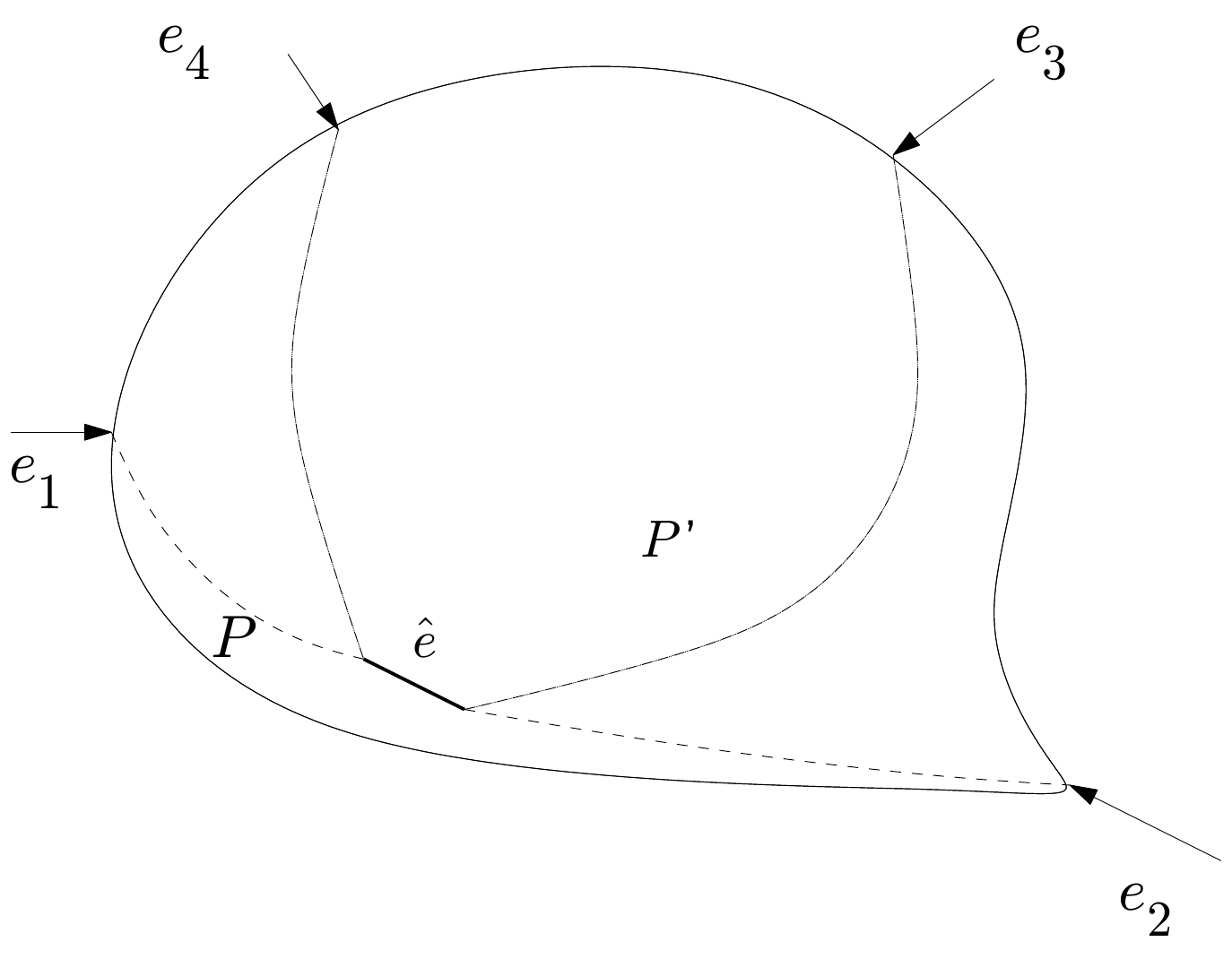}}
	\caption{The high edge $\hat e$ in the path $P$ 
joining $e_1$ and $e_2$}
\end{figure}

It follows that $C$ is incident to at least $4$ distinct half-edges of low edges of 
$2$-cells in $f^{-1}(\alpha_n)$, namely $e_1=\sigma(e_2)$, $e_2=\sigma(e_1)$, $e_3$ and $e_4$.
Suppose that these are the only $4$  half-edges of low edges of $2$-cells in $f^{-1}(\alpha_n)$ incident at $C$.  Then $\sigma(e_4)\in\{e_1,e_2,e_3\}$ and $\sigma(e_3)\in\{e_1,e_2,e_4\}$.   
There are essentially two cases.

\medskip\noindent{\bf Case 1.} 
Suppose first that $\sigma(e_4)=e_1$.  Then there is a path $Q'':=e_4\cdot P''\cdot e_1^{-1}$ with $P''$ in $C$, such that 
 $f(Q'')$  is an initial segment of the attaching path $R^{\pm n}$ of the $2$-cell $\a_n$.

Recall that $Q=e_1\cdot P\cdot e_2^{-1}$ and $Q':=e_3\cdot P'\cdot e_4^{-1}\ne Q$, where $P,P'$ are paths in $C$ that pass through the high edge $\hat e$.
Let $Q_1$ denote the part of $Q$ from  $e_1$ to the midpoint of $\hat{e}$, and $Q_2$ the part of $Q'$ from the midpoint of $\hat{e}$ to $e_4$.

\begin{figure}
	\centerline{\includegraphics[scale=0.3]{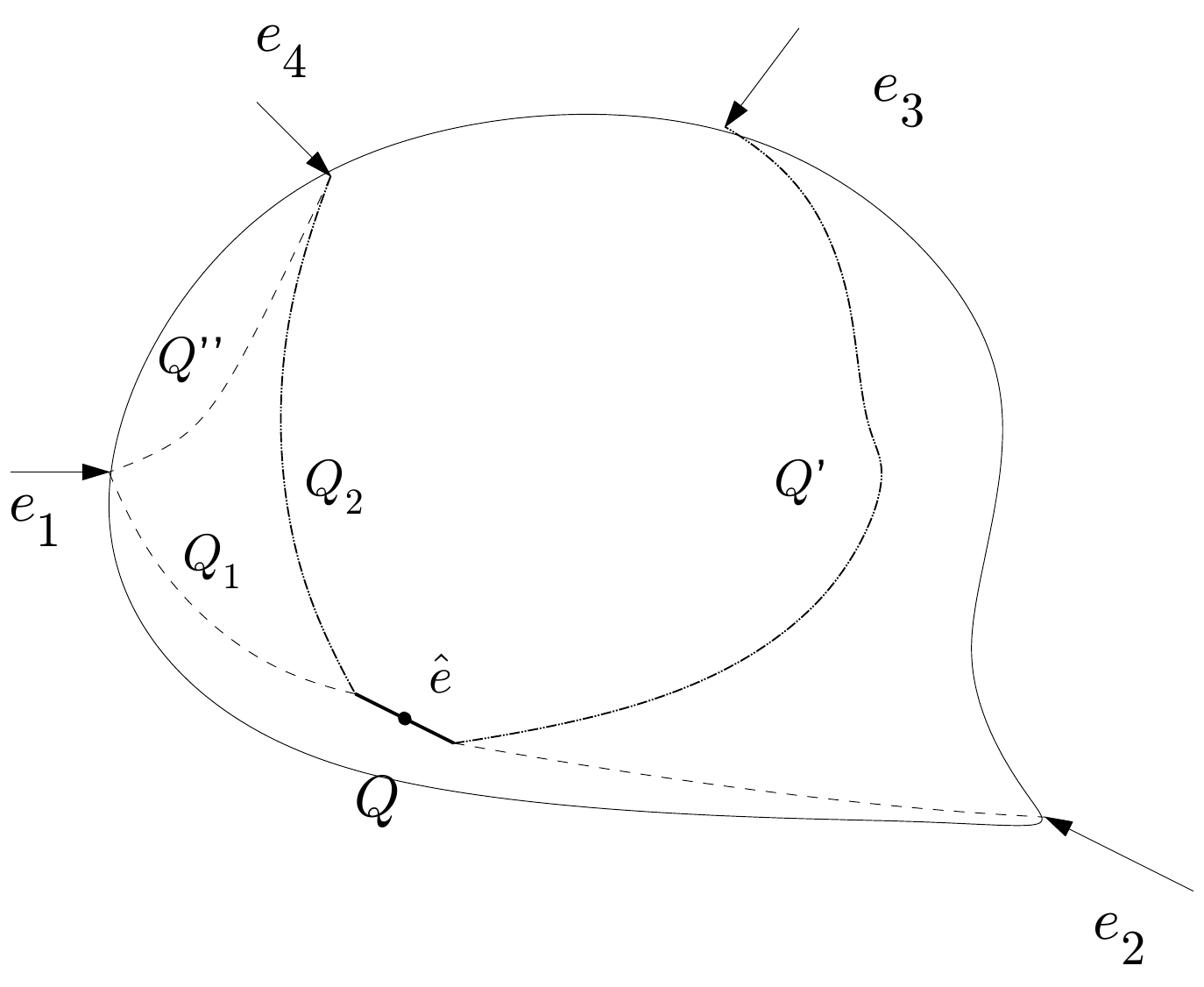}}
	\caption{Case 1 of a treeoid component of valence 4}
\end{figure}

Then $h_1:=(\psi_n\circ f)_*(Q_1)=U=V^{-1}$, $h_2:=(\psi_n\circ f)_*(Q_2)$ is a cyclic subword of $R^{\pm 1}$, and $h'':=(\psi_n\circ f)_*(Q'')$ is an initial segment of $(UV)^{\pm 1}$.  By  Lemma \ref{UV}  $h_1h_2$ and $h''$ are both positive, so $h_1h_2h''>1$, contradicting the hypothesis that  $(\psi_n\circ f)_*(\pi_1(C))=1$.

Similar arguments yield contradictions in the case where $\sigma(e_4)=e_2$, and in the cases where $\sigma(e_3)\in\{e_1,e_2\}$.

\medskip\noindent{\bf Case 2.} We are now reduced to the case where $\sigma(e_3)=e_4$ and $\sigma(e_4)=e_3$.  
As was the case for $Q$, $e_3$ and $e_4$ are low edges of the same $2$-cell.
and there is a path in $C$ joining them labelled $R^{\pm 1}$.
But the path $P'$ also joins these half-edges, and so the label on $Q'$ 
is $R^{\pm 1}$ and the $2$-cell concerned is $\alpha''$.

Arguing as in the case of $\hat{e}$, $Q'$ contains a high edge $\widehat {e'}$ of $\alpha''$. 
Since $\widehat{e'}$ is not a collapsible edge of $\alpha''$, it is contained in the attaching path of a $2$-cell of $f^{-1}(\alpha_n)$ that meets $C$ but not $e_3$ or $e_4$.  Thus $e_1$ or $e_2$ is a low edge of this other $2$-cell, and so this other $2$-cell must be $\alpha'$.   
Hence $\hat{e}$ is contained in $Q'$, while $\widehat{e'}$ is contained in $Q$.   
Let $L_1$ denote the subpath in $Q$ from the midpoint of $\hat{e}$ to $e_1$, $L_2$ the subpath of $Q$ from $e_1$ to the midpoint of $\widehat{e'}$, $L_3$ the subpath of $Q'$ from the midpoint of $\widehat{e'}$ to $e_3$, 
and $L_4$ the subpath of $Q'$ from $e_3$ to the midpoint of $\hat{e}$.  
Let $h_1,h_2,h_3,h_4$ denote the images of $L_1,L_2,L_3,L_4$ respectively in $G$ under $(\psi_n\circ f)_*$.   
Then $h_1=h_3=U^{-1}=V$, so  by Lemma \ref{UV} $h_1\le h_4^{-1}$ and $h_3\le h_2^{-1}$.  
Moreover these inequalities are strict because $\hat{e}\ne \widehat{e'}$ since these are high edges of distinct $2$-cells
 (as the edges are not almost-collapsible). 
Thus $h:=h_1^{-1}h_4^{-1}h_3^{-1}h_2^{-1}>1$, contradicting the hypothesis that $(\psi_n\circ f)_*(C)=1$.

\begin{figure}
	\centerline{\includegraphics[scale=0.3]{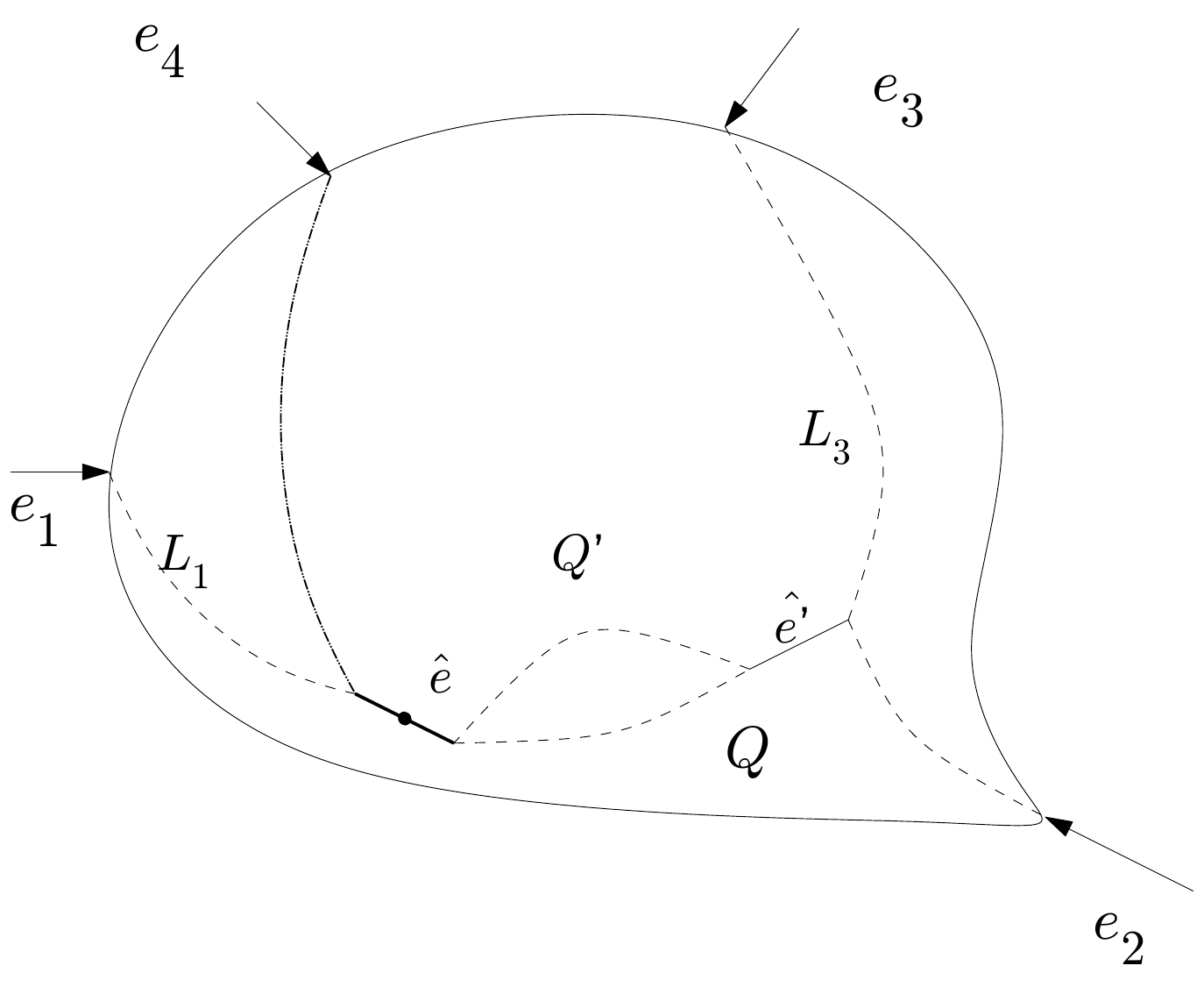}}
 	\caption{Case 2 of a treeoid component of valence 4}
\end{figure}

In all cases we have obtained a contradiction, and so $C$ meets at least $5$ half-edges of low edges of $2$-cells in $f^{-1}(\alpha_n)$, as claimed.

\medskip
To complete the proof, suppose that $Z''$ has $M_0$ treeoid components and $M_1$ non-treeoid components.   And suppose that there are $K$ $2$-cells in $f^{-1}(\alpha_n)$, and $J$ low edges of these $2$-cells.  Then $Y'$ is obtained from the $M_0+M_1$ components of $Z''$ by adding $J$ $1$-cells and $K$ $2$-cells.  The first Betti number of $Z''$ is at least $M_1$, since each non-treeoid component makes a positive contribution by definition.   Of the additional $1$-cells, $M_0+M_1-1$ are required to connect the different components of $Z''$, so the first Betti number of $Z''\cup (Y')^{(1)}$ is at least $J-M_0+1$.  Finally, each $2$-cell attached can reduce the first Betti number by at most $1$ so we have
\begin{equation}\label{e1}
\beta \ge J - M_0 + 1 - K.
\end{equation}

The $J$ low edges split into $2J$ half-edges, of which each treeoid component of $Z''$ meets at least $5$, as we have seen.  Hence $5M_0\le 2J$.  Combining this with inequality (\ref{e1}) and the fact that $J\ge 2K$, we have
$ 5\b > 5J -5K  - 5M_0 \ge 3J -5K\ge K$ as required.
\end{proof}
 
\begin{cor}\label{11k}
Let $g:Y'\imm \widehat{Y}_n$ be an immersion, where $n>1$. 
 Suppose that $Y'$ is compact and connected, 
 that no $2$-cell in $g^{-1}(\alpha_n)$ has an almost-collapsible edge, 
 and that $\pi_1(Y')$ can be generated by $k$ elements.   
 Then the number of $2$-cells in $g^{-1}(\alpha_n)$ is at most $11k$.
\end{cor}

\begin{proof}
Let $f:=\psi_n\circ g:Y'\to Y$.  Then $f$ is a branch map, and satisfies the conditions of Theorem \ref{notree}.  Let $p$ be a prime factor of $n>1$ and let $F$ be the field of order $p$.  Then the dimension over $F$ of $H_1(Y',F)\cong \pi_1(Y')^{ab}\otimes_\Z F$ is less than or equal to $k$, since $\pi_1(Y')$ can be generated by $k$ elements. By Theorem \ref{p-power}, at most $k$ of the $2$-cells in $f^{-1}(\alpha)$ are attached along $p$-th powers - in particular along $n$'th powers.

Now let $Y''$ be the subcomplex of $Y'$ obtained by removing all $2$-cells in $f^{-1}(\a)$ that are attached along $n$'th powers.  Now $Y'$ has first Betti number at most $k$ and is formed from $Y''$ by attaching at most $k$ $2$-cells, so $Y''$ has first Betti number at most $2k$.  Moreover, the restriction of $f$ to $Y''$ satisfies the hypotheses of Theorem \ref{5beta}, since $Y''$ contains no $2$-cell of $f^{-1}(\a)$ that is attached along an $n$-th power.  So the number of $2$-cells of $Y''\cap f^{-1}(\a)$ is at most $5$ times the first Betti number of $Y''$.   Hence the total number of $2$-cells in $f^{-1}(\a)$ is at most
$$5\times(2\times k) + k = 11k.$$
\end{proof}

\section*{Proof of Theorem A}

Following the proof of Theorem C, we reduce  to the two-factor case, where $|\Lambda|=2$.  
We restate and prove it in that form.

\begin{thm}[  = Theorem \ref{main}]
	Let $A,B$ be coherent locally indicable groups, and let $R=S^n\in A*B$, where $S$ is cyclically reduced of length at least $2$ and $n>1$.  Then $\G:=\frac{A*B}{\<\<R\>\>}$ is coherent.
\end{thm} 
 
The proof follows Wise's method in Section 4 of \cite{Wi}, beginning with a  lemma of Peter Scott:
\begin{lem}[Scott]\label{Peter}
	Let $H$ be a freely indecomposable group generated by $n$ elements.  
	Assume that every subgroup of $H$ generated by fewer than $n$ elements is finitely presentable.  
	Then there exists a finitely presentable, freely indecomposable  group $K$ and an epimorphism $K\onto H$ that does not factor through a proper free product: if $K\onto H'\onto H$ then $H'$ is freely indecomposable.
\end{lem}

\noindent{\bf Remark.}  Scott \cite[Lemma 2.2]{Sc} proves this in 
the context of his inductive proof of coherence of $3$-manifold groups.  
The inductive assumption on subgroups of $H$ is used in an essential way in the proof of 
Lemma \ref{Peter}, so to avoid logical difficulties we shall be careful to use it also in the context of an inductive proof of Theorem \ref{main}.  
 While preparing this work for print, Louder and Wilton brought to our attention Swarup's preprint \cite{Sw} 
	where he gives a proof he attributes to Delzant of the general case of this result (without 
	the condition that lower rank subgroups be finitely presented).

\bigskip
Let $H$ be a finitely generated subgroup of $\G$.
We want to show that $H$ is finitely presentable.
We proceed by induction on the number of generators for $H$:
clearly it is true for subgroups with at most 1 generator. 
We can assume that $H$ is generated by $k$ elements, and that every subgroup of $G$ generated by strictly fewer than $k$ elements is finitely presentable.  In particular, if $H$ decomposes as a free product, then each free factor (and hence also $H$) is finitely presentable, using Grushko's theorem.  So we may assume that $H$ is freely indecomposable.

First construct a $2$-complex $\widehat{Y}$ (resp $Y$) with fundamental group $\G$ (resp. $G:=(A*B)/\<\<S\>\>$) as follows.  Let $X_A,X_B$ be presentation $2$-complexes for $A,B$ respectively.  
Add a $1$-cell $e$ to $X:=X_A\sqcup X_B$ joining the base-points of the two components, and then attach a $2$-cell $\alpha$ along $\sigma^n$ (resp. $\sigma$), 
where $\sigma$ is a path in $X^{(1)}\cup e$ representing $S\in A*B = \pi_1(X\cup e$).  Thus $Y$ is a simple enlargement of $X$, and $\widehat{Y}\to Y$ is the $n$-fold branched cover of $Y$, as defined in \S \ref{techstuff}.

\smallskip
\noindent{\bf Definition.}
Say that an edge in a  2-complex is  {\em \inac} if it occurs in the boundary of
at least one 2-cell, and is not almost-collapsible, in particular not collapsible. 
This means that either the edge occurs at least once each in the boundaries of two different 2-cells $\alpha_1, \alpha_2$
i.e. such that reading the boundaries $\partial \alpha_1,\partial\alpha_2$ from the edge $e$ 
with the orientation induced by a choice of orientation for $e$ gives two different  words,
or the edge occurs   twice on the boundary of a single 2-cell $\alpha_1$, such that reading 
the boundary $\partial\alpha_1$ starting from the different occurrences of the edge $e$ 
with the orientation induced by a choice of orientation for $e$ gives two different words.

\begin{lem}
	Let $H$ be  a $k$-generated freely indecomposable subgroup of $G$.

	There is a finitely presented, $k$-generated group $K$, 
	and a finite 2-complex $Y_0$ with $\pi_1(Y_0)=K$,
	where all edges of $Y_0$ are \inac,  $Y_0 \imm Y$,
	and the immersion induces  an epimorphism $F:\pi_1(Y_0)\onto H$.
\end{lem}

\begin{proof}
Let $K$ and   
$ \f: K \to H$ be the finitely presented group and epimorphism 
provided by  Scott's Lemma \ref{Peter}.
Let $W_0$ be a finite 2-complex with $\pi_1(W_0) =K$, let $q_H:Y_H\to \widehat{Y}$ be the covering of $\widehat{Y}$ with $\pi_1(Y_H)=H$,
and subdivide the $1$- and $2$-cells  so that 
the map $\f:K\to H$ is induced by a combinatorial map $g_0: W_0 \to Y_H$. 
Now $g_0$ is not in general an immersion, but by Lemma \ref{fold} it can be factored as $W_0\to Y' \imm Y_H$ 
where $W_0\to Y'$ is surjective
and  $\pi_1$-surjective, $\pi_1(Y')$ is indecomposable, and $f_0:Y'\imm Y_H$ is an immersion.  

Thus there exists at least one $\pi_1$-surjective immersion $Y'\imm Y_H$ 
such that $Y'$ is compact and connected with $k$-generated, freely indecomposable fundamental group.  
Among all such, choose one $f_0:Y_0\imm Y_H$ with $Y_0$ minimal, in the sense of having fewest possible cells.  
(Note that this choice may entail a change of basepoint in $Y_H$, thus effectively replacing $H$ by a conjugate.)

\noindent 
If $Y_0$ had a non \inac\ edge $e$,  then either 
\begin{enumerate}
	\item $e$ does not lie in the boundary of any 2-cell and
	\begin{enumerate}
		\item $e$	separates,  giving either 
		a free decomposition of $\pi_1(Y_0)$, or a smaller subcomplex carrying $K$;
		\item $e$ does not separate, and there is a free $\Z$ factor.
	\end{enumerate}
	\item $e$ is collapsible, and there is a  smaller subcomplex carrying $K$;
	\item $e$ is almost-collapsible and $K$ has a finite cyclic factor 
\end{enumerate}

Conclusion: every edge of $Y_0$ is \inac.  
\end{proof}

\begin{lem}\label{infseqimm}
	Let $Y=X\cup e\cup\a$ be a simple enlargement of a $2$-complex $X$ with locally indicable fundamental groups, and $\wh{Y}$  its $n$-fold branched cover, for some $n\ge 2$.  Let $H$ be a freely indecomposable, finitely generated subgroup of $\pi_1(\wh{Y})$, and let $Y_H$ be the connected cover of $\wh{Y}$ with fundamental group $H$.  Then there exists a sequence of $\pi_1$-surjective immersions
	\begin{equation}\label{imms}
		Y_0 \looparrowright Y_{1}  \looparrowright Y_{2} \dots Y_i \looparrowright Y_{i+1} \looparrowright \dots   \to Y_H
	\end{equation}
	such that 
	\begin{enumerate}
		\item[(a)] each $Y_j$ is a finite connected $2$-complex;
		\item[(b)]  every $1$-cell of $Y_j$ is \inac;
		\item[(c)] Each $Y_{j+1}$ is the union of the image of $Y_j$ and of a  reduced 
van Kampen diagram $\delta_j:D^2\to Y_{j+1}$ with $\delta_j(\partial D^2)$ in the image of $Y_0$;
		\item[(d)] each immersion $Y_j\imm Y_H$ is also $\pi_1$-surjective; and 
		\item[(e)] each element of $\mathrm{Ker}(\pi_1(Y_0)\to\pi_1(Y_H))$ is mapped to the identity element of $\pi_1(Y_j)$ for some $j$.
	\end{enumerate}
\end{lem}

\begin{proof}
Note first that, for any $2$-complex $Z$ and immersion $f:Z\imm Y_H$, if $e$ is a \inac\ edge of $Z$ then its image $f(e)$ is \inac\ in the image $f(Z)\subset Y_H$.  Hence if $f$ factors through a surjective map $Z\twoheadrightarrow Z'$, then the image of $e$ in $Z'$ will be \inac\ in $Z'$.   Note also that, if $\delta:\Delta\to Y_H$ is a  reduced  
van Kampen diagram, then any interior edge of $\Delta$ is \inac\ in $\Delta$. (Otherwise the two  $2$-cells incident to the edge   would cancel, contradicting the hypothesis that the diagram is reduced.)

We shall create the complexes $Y_j$, beginning with the complex $Y_0$ given by Lemma 5.6, in such a way that each edge is the image either of an edge of $Y_0$ or of an interior edge of a   reduced
van Kampen diagram.  So the fact that their edges are all \inac\ follows from the above remarks.

Let $f_0:Y_0\imm Y_H$ be the immersion given by Lemma 5.6.	
		Choose a sequence of  cyclically reduced closed paths $w_0, w_1, \dots $ in $Y_0^{(1)}$ which together normally generate the kernel of  $ f_{0,*}:\pi_1(Y_0)\to\pi_1(Y_H)$, 
		and a  reduced 
van Kampen diagram $\phi_j:D_j\to Y_H$ for each $w_j=1$ over $\frac{A*B}{\<\<R\>\>}$
	(we have fixed presentations for $A$ and $B$).

	Form the adjunction space  $W_1=Y_0\cup\partial D_0$ by identifying 
	$\delta_0=\partial D_0$ with its image in $Y_0$, 
	so that the maps $f_0,\phi_0$ combine to give a map $W_1\to Y_H$.  

Then apply Lemma \ref{fold} to factor this map as the composite of a surjective, $\pi_1$-surjective map $W_1\to Y_1$ and an immersion $f_1:Y_1\imm Y_H$.  By our initial remarks we see that every edge of $W_1$ and of $Y_1$ is \inac\ in $Y_1$.
	Finally note that $w_0\in\mathrm{Ker}(\pi_1(Y_0)\to\pi_1(Y_1))$.

Now iterate this construction.  Suppose inductively that we have defined $Y_i$, 
	and immersions $Y_0\looparrowright \cdots\looparrowright  Y_i$, such that every edge of each of $Y_0,\dots,Y_i$ is \inac\ and
	$w_0,\dots,w_{i-1}\in\mathrm{Ker}(\pi_1(Y_0)\to\pi_1(Y_i))$.  Then we can adjoin the
	diagram $D_i$ to $Y_i$ and factor the resulting map to $Y_H$ through an immersion, giving a new 2-complex $Y_{i+1}$ and an immersion $Y_i \looparrowright Y_{i+1}$
	that factorises $f_i:Y_i\looparrowright Y_H$ through $f_{i+1}:Y_{i+1}\looparrowright Y_H$\,, and $w_i\in\mathrm{Ker}(\pi_1(Y_0)\to\pi_1(Y_{i+1}))$.   
	 Again, arguing as above, every edge in $Y_{i+1}$ is \inac\ in $Y_{i+1}$.

In this sequence of immersions, each $f_j:Y_j\imm Y_H$ is $\pi_1$-surjective, since the $\pi_1$-surjective immersion $f_0$ factors through $f_j$.
Moreover, any element of $\mathrm{Ker}(f_{0,*})$ lies in the normal closure of the $w_i$, and hence is mapped to the identity element of $\pi_1(Y_j)$ for some $j$, as required.   
\end{proof}

\noindent{\bf Definition.}  The {\em $R$-subcomplex} $R(Y_j)$ of $Y_j$ is the union of the image of $Y_0$ in $Y_j$ together with the preimages of $\a$ under $Y_j\imm Y$ and all of their incident $0$-and $1$-cells.
If $\delta:\Delta\to Y$ is a van Kampen diagram, then the {\em $R$-subcomplex} $R(\Delta)$ of $\Delta$ is the union of $\partial\Delta$ with the preimages of $\a$ under $\delta:\Delta\to Y$ and all of their incident $0$-and $1$-cells.

\begin{lem}\label{Rconn}
	The complexes $R(Y_j)$ are all connected.
\end{lem}

\begin{proof} This follows by induction, starting from the fact that  $Y_0$ is connected.
	
	The  complex $Y_{j+1}$ is the image of $W_{j+1}$, which itself is formed from $Y_j$ 
	by gluing on the   reduced 
van Kampen diagram $D_j$.
	The fact that $R(D_j)$ is connected for a  reduced 
van Kampen diagram for $w_j$
	is one  of the essential properties established in the proof
	of the Freiheitssatz (Theorem \ref{Freiheitssatz})  for one relator products of locally indicable groups in \cite{Br0,H1,Sh},
	that the natural map $A \to \frac{A*B}{\<\<R\>\>}$ is an injection (and thus
	that for each $i$ the word $w_i$ labelling  the boundary $\delta_i$ 
	contains non--trivial words in $A$ and in $B$).  
	If $R(D_j)$ were not connected, then there would be a subdiagram of $D_j$
	containing some regions labelled $R$, and whose boundary is a word in $A$ or in $B$. 
  But in a reduced diagram, any subdiagram with boundary label in $A$ or $B$ cannot contain $R$-regions, so this is a contradiction.
	
	Thus $R(Y_{j+1})$ is the image of the union of two connected subcomplexes of $Y_j$ with non empty intersection.
\end{proof}

\begin{lem}\label{subseq}
	There is an infinite subsequence
	\begin{equation}\label{Ysigma}
Y_{\sigma(0)}\imm Y_{\sigma(1)}\imm \cdots Y_{\sigma(j)}\imm\cdots\imm Y_H
\end{equation}
	of (\ref{imms}) of Lemma \ref{infseqimm} such that the restrictions
	\begin{equation}\label{Rsigma}
R(Y_{\sigma(0)})\imm R(Y_{\sigma(1)})\imm \cdots R(Y_{\sigma(j)})\imm\cdots
\end{equation}
	are all cellular isomorphisms.
\end{lem}

\begin{proof}
Let $Z_i =R(Y_i)$;
it is clear that the map  $Z_i\imm Y_H$ is still $\pi_1$-surjective 
(as the initial map $f_0:Y_0 \imm Y_H$ factors as $Y_0 \imm Z_i \imm Y_H$).

By Corollary \ref{11k} (every edge in $Y_i$ is \inac) applied to the composite $q_H\circ f_i:Y_i\looparrowright Y_H \to Y$, the number of preimages of $\a$ in $Y_i$ is bounded above by $11k$, where $k$ is the minimal number of generators of $H$. Since every $1$-cell in $Y_i$ is incident at a $2$-cell, the number of preimages of $e$ in $Y_i$ is also bounded (for example, 
by $11k\ell$ where $\ell$ is the length of $R$ as a path in $Y^{(1)}$).

But every cell in $Z_i$ is either in the image of the compact $2$-complex $Y_0$, or is incident to one of the boundedly many $2$-cells that are preimages of $\a$. 
Hence
the $Z_i$ are bounded in size, so some subsequence $\{Z_{\sigma(i)}\}$ of the $Z_i$ stabilises,
in the sense that  the sequence $Z_{\sigma(1)} \imm Z_{\sigma(2)} \imm ...$ consists of isomorphic $2$-complexes.   
Now Wise \cite[Lemma 4.2]{Wi} has shown that any immersion from a finite complex to itself must be an isomorphism.  
(Otherwise some power of the immersion is a retraction onto a proper subcomplex, 
contrary to the definition of immersion; we are grateful to Lars Louder  for explaining this to us.). 
Hence in fact we may assume that
$Z_{\sigma(1)} \imm Z_{\sigma(2)} \imm ...$  consists of isomorphisms, as claimed.  
\end{proof}

\noindent{\bf Last step  of the proof of Theorem A} 

By Lemma \ref{subseq} there is a subsequence (\ref{Ysigma}) of (\ref{imms}) such that all the maps $R(Y_{\sigma(i)})\imm R(Y_{\sigma(i+1)})$ in (\ref{Rsigma}) are cellular isomorphisms.
Using this, we set $Z:=R(Y_{\sigma(1)})$ and identify each $R(Y_{\sigma(i)})$ with $Z$ via the inverse isomorphism.

Any element of the kernel of the map $\pi_1(Z)\to \pi_1(Y_H)=H$ becomes trivial in some $\pi_1(Y_{\sigma(i)})$.   
But $Z=R(Y_{\sigma(i)})$ is a subcomplex of $Y_{\sigma(i)}$, 
so any such element can be expressed as the boundary label of a van Kampen diagram $\Delta$ in $Y_{\sigma(i)}$.  
As in Lemma \ref{Rconn}, the $R$-subcomplex of $\Delta$ may be assumed to be connected, by Theorem \ref{Freiheitssatz}.  Thus its complement consists of a finite number of open disks in $\Delta$ whose boundaries are paths in $Z$ and whose images under the map $(q_H\circ f_{\sigma(1)}):Y_{\sigma(1)}\to Y$ lie in $X_A$ or in $X_B$.   Since $Z\cap (q_H\circ f_{\sigma(1)})^{-1}(X_A\sqcup X_B)$ is compact, and each of $A,B$ is coherent, it follows that $Ker(\pi_1(Z)\to H)$ is normally generated by a finite number of loops in $Z\cap (q_H\circ f_{\sigma(1)})^{-1}(X_A\sqcup X_B)$.  
For sufficiently large $N$, these loops are all nullhomotopic in $Y_{\sigma(N)}$, and it follows that $\pi_1(Y_{\sigma(N)})\to H$ is an isomorphism.  
Since $Y_{\sigma(N)}$ is compact, $H$ is finitely presentable, as claimed.
\qed

\smallskip
We end this section by noting that stronger versions of 
 some of our results -- notably Theorem 
\ref{BrodskiiLemma} -- hold in the torsion case.  
Translated into our set-up, a result in \cite{DH} says the following.

\begin{thm}[\cite{DH}, Theorem 3.3]\label{DHthm}
Let $Y:=X\cup e\cup\alpha$ be a simple enlargement where $X$ is a $2$-complex, 
all of whose fundamental groups are locally indicable.
Let $n>1$ be an integer, and let $\widehat{Y}_n\to Y$ be the $n$-fold cyclic branched cover defined at the start of this section.  Let $\mathcal{P}:\Sigma\to \widehat{Y}_n$ be a reduced picture on a compact orientable surface $\Sigma$ that contains $V$ $\alpha$-discs.  Then the number of points where $\mathcal P$ meets $\partial\Sigma$ is at least
 $$2V(n-1) +2\chi(\Sigma) \ .$$
\end{thm}

\begin{cor}\label{torsionBrodskii}
Let $$G:=\frac{*_{\lambda\in\Lambda} G_\lambda}{\<\<R^n\>\>},$$
where $G_\lambda$ are locally indicable groups, $R\in *_\lambda G_\lambda$ is cyclically reduced of length at least $2$, and $n>1$ is an integer.  Then
\begin{enumerate}
\item If $\lambda,\mu\in\Lambda$ and $g\in G$ are such that the intersection $G_\lambda\cap G_\mu^g$ in $G$ is non-trivial, then $\mu=\lambda$ and $g\in G_\lambda$.
\item Any free abelian subgroup of $G$ of rank greater than $1$ is contained in a conjugate of $G_\lambda$ for some $\lambda\in\Lambda$.
\item Any subgroup of $G$ that is isomorphic to a   right angled Artin group based on a connected graph with at least one edge is contained in a conjugate of $G_\lambda$ for some $\lambda\in\Lambda$.
\end{enumerate}
\end{cor} 

\begin{proof}
As usual (see proof of Theorem C) we work in the two-factor case  $(A*B)/\<\<R^n\>\>$,
with  $A:=*_{\lambda\in\Lambda(1)} G_\lambda$,   $B:=*_{\lambda\in\Lambda(2)} G_\lambda$.
Without loss of generality we can assume that $R$ is not a proper power in $A*B=*_\lambda G_\lambda$.

We then form a simple enlargement $Y=X\cup e\cup\alpha$ of $X:=X_A\sqcup X_B$ where $X_A,X_B$ are connected $2$-complexes; $\pi_1(X_A)\cong A$ and $\pi_1(X_B)\cong B$,
 $e$ is a $1$-cell joining the base points of $X_A$ and $X_B$, and $\alpha$ is a $2$-cell attached along a path representing $R\in A*B$.  
The result follows by applying Theorem \ref{DHthm} to pictures over $\widehat{Y}_n$ as follows.

\begin{enumerate}
\item A conjugacy relation $y=x^g$ with $x\in G_\lambda$ and $y\in G_\mu$ can be expressed using a reduced picture $\mathcal{P}:\Sigma\to\widehat{Y}_n$ where $\Sigma$ is an annulus and the components of $\partial\Sigma$ map to reduced paths in $X$ which represent $x,y$ respectively.  In particular these paths do not involve $e$, and so no arc of $\mathcal{P}$ meets $\partial\Sigma$.  Since $\chi(\Sigma)=0$, it follows from Theorem \ref{DHthm} that $\mathcal{P}$ has no $\alpha$-discs.  So we can regard $\mathcal{P}$ as a picture over $X\cup e$, and the conjugacy relation $y=x^g$ already holds in the free product $\pi_1(X\cup e)=A*B=*_\lambda G_\lambda$.  The result follows from well-known properties of conjugacy in free products.
\item A commutator relation $xy=yx$ in $G$ can be expressed using a reduced picture $\mathcal{P}:\Sigma\to\widehat{Y}$ where $\Sigma$ is the torus.   
Since $\chi(\Sigma)=0$ and $\partial\Sigma=\emptyset$, Theorem \ref{DHthm} again tells us that $\widehat{P}$ has no $\alpha$-discs.  Thus $\widehat{P}$ is a picture over $X\cup e$, and the commutator relation $xy=yx$ already holds in $\pi_1(X\cup e)=A*B=*_\lambda G_\lambda$.   
The result again follows from well-known properties of commutation in free products. 

\item The final assertion follows easily from the first two.  The subgroup corresponding to any edge is contained in some $G_\lambda^g$, and the subgroup corresponding to any vertex can be contained in at most one $G_\lambda^g$, so there is a unique $\lambda\in\Lambda$ and a unique right coset $G_\lambda g$ such that $G_\lambda^g$ contains the right angled Artin group in question.
\end{enumerate}
\end{proof}

\section{Wreath products and Thompson's group as subgroups}\label{Z3}

In this section we prove the first  two parts of Theorem  \ref{noThompson}.

\begin{lem}\label{refslemma}
Let $G$ be as in the Theorem, and let $K:=\<H,x\><G$, where $H<G_\lambda^g$ for some $\lambda\in\Lambda$ and $g\in G$.  If $H^x\cap H$ is non-cyclic, then $K<G_\lambda^g$.
\end{lem}

\begin{proof}
$H^x\cap H<G_\lambda^{gx}\cap G_\lambda^g$, so $x\in G_\lambda^g$ by Theorem \ref{BrodskiiLemma}.
\end{proof}

\begin{pfof}{Theorem \ref{noThompson} (1) and (2)}

\medskip\noindent{\bf Case  1.} $K\cong \Z~wr~\Z$.

The commutator subgroup $[K,K]$ of $K$ is free abelian of infinite rank, 
with basis $\{x_n,~n\in\Z\}$, and $K=[K,K]\rtimes\<t\>$ 
where $tx_n=x_{n+1}t$ for all $n$.   Define $H:=\<x_1,x_2,x_3\><K$, and note that $K=\<H,t\>$ and $H\cap H^t=\<x_1,x_2\>$ is non-cyclic.  By Corollary \ref{noZ3} there is a unique $\lambda\in\Lambda$ and a unique coset $G_\lambda g$ such that $H<G_\lambda^g$.  Applying Lemma \ref{refslemma} with $x=t$  gives the result.  

\medskip\noindent{\bf Case  2.} $K$ is isomorphic to Thompson's group $F$.

There is a presentation of $K$ of the form
$$K=\<~x_1,x_2,x_3,\dots~|~x_mx_n=x_{n+1}x_m~(m<n)~\>.$$
Let $y_n:=x_n^{-1}x_{n+1}$.  Then $y_m$ commutes with $x_n$ and $y_n$ in $K$ if $m+1<n$.  
Taking $H=\<y_k,y_{k+2},y_{k+4}\>$ for any value of $k$, we see that $H<G_\lambda^g$ for a unique $\lambda\in\Lambda$ and a unique coset $G_\lambda g$.   By Lemma \ref{refslemma} it follows that $y_\ell\in G_\lambda^g$ for every $\ell\ge k+6$ and every $0<\ell\le k-2$.  Applying  these facts for varying $k$ gives $y_j\in G_\lambda^g$ for all $j$, and another application of Lemma \ref{refslemma} then gives $x_1\in G_\lambda^g$ and hence $K<G_\lambda^g$ as required.
\end{pfof}

\section{ A new proof of  {\Brodsky}'s lemma} \label{Blemma}

\begin{pfof}{Theorem \ref{BrodskiiLemma}}
By Theorem \ref{torsionBrodskii} we may assume that $G$ is torsion-free, that is that the relator $R$ is not a proper power in $*_\lambda G_\lambda$.
By Corollary \ref{Shapiro} we know that the intersection $G_\lambda\cap G_\mu^g$ is free unless the conclusion of the theorem holds.  So let us suppose that there is a free subgroup  $K$ of rank $2$ in $G_\lambda\cap G_\mu^g$.

 Now put $A:=G_\lambda$ and $B:=\ast_{\nu\ne\lambda} G_\nu$, 
so that $G$ is a one-relator product of the two locally indicable groups $A,B$. 
 Then we model the situation geometrically as :
  $Y=X\cup e\cup\alpha$ is a simple enlargement of $X:=X_A\sqcup X_B$ where $X_A$ and $X_B$ are connected $2$-complexes with locally indicable fundamental groups $A,B$ respectively and $e$ is a $1$-cell and $\alpha$ is a $2$-cell attached along a path representing $R\in \pi_1(X\cup e)\cong A*B\cong *_\lambda G_\lambda$.

 Now $K$ is a free subgroup of $A$ of rank $2$. 
 Suppose that  $P$ is either a loop at the base point of $X_A$ or a path from the base point of $X_A$ to the base-point of $X_B$, such that $P^{-1}KP$ is equal in $\pi_1(Y) $  to a subgroup of $A$ or of $B$
respectively.   Let $V=S^1\vee S^1$ be a $2$-petal rose with base-point $*$, and let $\phi:V\times [0,1]\to Y$ be a map representing this conjugacy set-up. That is, $\phi(V\times\{0\})$ is a pair of curves in $X_A$ that generate $K$, $\phi(V\times\{1\})$ is a pair of curves in $X$, and $\phi(\{*\}\times [0,1])$ is the path $P$.

After subdividing cells in $V\times [0,1]$ to make $\phi$ be cellular, we can apply Lemma \ref{fold} to factor $\phi$ as a surjective, $\pi_1$-surjective map $\phi':V\times [0,1]\to Y'$ followed by an immersion $f:Y'\imm Y$.  Let $\alpha$ denote the $2$-cell in $Y\setminus X$, $Z'$ the subcomplex of $Y'$ obtained by deleting every $2$-cell in $f^{-1}(\alpha)$ along with its low edge, and $X':=f^{-1}(X)\subset Z'$.

Now as well as being $\pi_1$-surjective, 
$\phi':V\times [0,1]\to Y'$ is  also $\pi_1$-injective since the composite $\phi=f\circ\phi':V\times [0,1]\to Y$ is $\pi_1$-injective.   Hence $\pi_1(Y')\cong\pi_1(V\times [0,1])$ is free of rank $2$.

The inclusion $V\times\{0\}$ into $V\times [0,1]$ is a homotopy equivalence, and $\phi'(V\times\{0\})$ is contained in a component $X_0$ of $X'$.  Since $V\times\{0\}\to Y'$ is $\pi_1$-surjective and factors through $X_0$, the inclusion $X_0\to Y'$ is also $\pi_1$-surjective.  Thus $X_0$ has first Betti number $\ge 2$.  The same applies to the component $X_1$ of $X'$ that contains $\phi'(V\times\{1\})$.

No component of $Z'$ has first Betti number $0$.  
For otherwise $Y'$ collapses onto that component, by Theorem \ref{notreetf}, 
and has first Betti number $0$, a contradiction.  
Since $Y'$ is connected with first Betti number $2$, and is formed from $Z'$ by attaching equal numbers of $1$- and $2$-cells, it follows that all but at most one of the components of $Z'$ have first Betti number $1$, and no component has first Betti number greater than $2$.

Since $Z'$ is formed from $X'$ by attaching $1$-cells, it follows in turn that at most one component of $X'$ can have first Betti number greater than $1$.  Hence $X_1=X_0$.

 $P':=\phi'(\{*\}\times [0,1])$ is a path in $Y'$ connecting two points of $X_0$, and $P=f(P')$, 
 so both ends of $P$ lie in $X_A$.  Thus $P$ is a loop in $Y$ based at the base-point of $X_A$.  
 Moreover, since $X_0$ is connected and $X_0\to Y'$ is $\pi_1$-surjective, 
 there is a path $P_0$ in $X_0$ with the same endpoints as $P'$, 
 and a loop $P''$ in $X_0$ that is homotopic rel base point in $Y'$ to $P'\cdot  P_0^{-1}$.  
 Thus the element of $G$ represented by $P$ is 
$$g:=[P]=[f(P')]=[f(P''\cdot P_0)]\in \pi_1(X_A)=A.$$ 
\end{pfof}

\section{Direct products of free groups as subgroups}\label{FxF}

\begin{pfof}{Theorem \ref{noThompson} (3)}
 Since $F_2\times F_2$ is freely indecomposable, we may reduce to the case where $|\Lambda|=2$, as in the proof of  Theorem \ref{sb2}.  As usual we model the situation using a
simple enlargement 
	$Y=X\cup e\cup\alpha$  of  $X:=X_A\sqcup X_B$ where $X_A$ and $X_B$ are connected $2$-complexes with locally indicable fundamental groups $A,B$ respectably.

Suppose that $a,b,c,d\in G$ generate $K$, a direct product  
of two free groups $\<a,b\>$ and $\<c,d\>$ of rank $2$. 
We shall show that $K$ is contained in  a conjugate of one of the factors $A$ or $B$ 
in $G$.

Let $V=S^1\vee S^1$ be a $2$-petal rose with base-point $*$.  
There is a map $\phi:V\times V\to  Y$, combinatorial after subdivision, 
 with $\phi_*(\pi_1(V\times V))= K < G$.  
The restriction to any one of the four subcomplexes of $V\times V$ 
isomorphic to $S^1\times S^1$ is a toral picture over $Y$ representing one of the 
four commutator equations $ac=ca$ etc.  
We may assume without loss of generality that each of the four pictures is reduced.

If $R$ is a proper power, then by  Theorem \ref{DHthm} 
none of the toral pictures contains a vertex.  
In that case $\phi(V\times V)\subset X\cup e$, so the embedding 
$\phi_*:\pi_1(V\times V)\to \pi_1(Y)$ factors through $\pi_1(X \cup e)=A*B$.  
Since $F_2\times F_2$ is freely indecomposable, $\phi_*(\pi_1(V\times V))$ must belong to a conjugate of $A$ or $B$.

Hence we may assume that $R$ is not a proper power.  
Using Lemma \ref{fold}, we factor $\phi$ through an immersion $f:Y'\imm Y$ 
such that $V\times V\to Y'$ is surjective and $\pi_1$-surjective.  
In particular $\pi_1(Y')\cong F_2\times F_2$ and $f:Y'\to Y$ is $\pi_1$-injective, 
since $\phi$ is $\pi_1$-injective.  
Thus $Y'$ is a connected $2$-complex.
  Moreover $H_1(Y')\cong H_1(F_2\times F_2)\cong\Z^4$ and $H_2(Y')$ surjects onto $H_2(F_2\times F_2)\cong\Z^4$,
so $\chi(Y')>0$.  Let $Z'\subset Y'$ be the subcomplex in Theorem \ref{notreetf}. 

Then $Z'$ is obtained from $Y'$ by the removal of equal numbers of $1$- and $2$-cells, 
so $\chi(Z')=\chi(Y')>0$,
and some component $T$ of $Z'$ has positive Euler characteristic.   
By Theorem \ref{notreetf} we may assume that $f_*(\pi_1(T))\neq\{1\}$ in $G$, 
for otherwise $Y'$ collapses to $T$ implying that  $f_*(\pi_1(Y))=f_*(\pi_1(T))=\{1\}$ in $G$, contrary to hypothesis.

Now $f(T)\subset X\cup e$. Since $\pi_1(X\cup e)=A*B$, 
the image $K$ of $\pi_1(T)$ in $\pi_1(X\cup e)$ splits as a free product $\ast K_j$, 
where each $K_j$ is either cyclic or contained in a conjugate of $A$ or $B$.
Suppose that some $K_j$ is contained in $A^g$. 
But as the factor groups are locally indicable, $A$ embeds in $G$ via the natural map
by Theorem \ref{Freiheitssatz}, 
and hence $K_j$ is isomorphic to a finitely generated subgroup of $F_2\times F_2$.  
Moreover, it follows from  
Corollary \ref{pi1-inj}  (1)  
that  $f|_{T}$ is $\pi_1$-injective,
and $\pi_1(T)\cong K_j$ is isomorphic to a finitely presented subgroup of $F_2\times F_2$
contained in $A^g$.

Observe that a finitely generated  subgroup $L$ of $F_2\times F_2$ 
either   contains a  copy of $F_2\times F_2$ or has Euler characteristic $\le 0$
(and is $\Z$,   $\Z\times\Z$ or free or  $\Z\times free$).
This follows easily by considering the restrictions to $L$ of the coordinate projections
$p_1,p_2:  F_2 \times F_2 \to F_2$. 
If either is injective, then $L$ is free.  
If  $p_1$ has rank 1 kernel, then there is 
an exact sequence  $1 \to \Z \to L \to p_1(L) \to 1$,
which splits since $p_1(L)$ is free.  
So $L \cong \Z \times free$.  

Finally, if $\ker(p_1)$ and $\ker(p_2)$ are both (free) of  rank $\ge 2$, 
then $L \supset \ker(p_1) \times \ker(p_2) \supset F_2 \times F_2$.

The components of $Z'$ with fundamental groups that are
free products of free groups and $\Z\times free$ groups  have
non-positive Euler Characteristic.
But we have seen that there is at least one component with positive 
Euler characteristic, so  some $\pi_1(T)$, and so $A^g$ 
contains a subgroup of $K$ isomorphic to $F\times F$.

Now any subgroup of $K=\<a,b\>\times\<c,d\>$ that is  isomorphic to $F_2\times F_2$ has the form $L=L_1\times L_2$ where $L_1,L_2$ are subgroups of $\<a,b\>,\<c,d\>$ respectively.   If $L<A^g$ then we have $F_2\cong L_2<A^g\cap A^{ga}$ and so $a\in A^g$ by \Brodsky's Lemma, Theorem \ref{BrodskiiLemma}.   Similarly $b,c,d\in A^g$ and so $K<A^g$. 
 Since $A$ is a free product of some of the $G_\lambda$, and $K$ is freely decomposable, it follows that 
$K<G_\lambda^g$ for some $g,\lambda$, as claimed.
\end{pfof}

\end{document}